\documentclass[12pt]{article}
\usepackage{amssymb,amsmath,eucal}

\renewcommand{\a}{\alpha}
\renewcommand{\b}{\beta}

\newcommand{\e}{\epsilon}
\newcommand{\f}{\varphi}
\newcommand{\g}{\gamma}

\renewcommand{\i}{\iota}
\renewcommand{\j}{\psi}

\renewcommand{\l}{\lambda}
\renewcommand{\o}{\omega}

\newcommand{\s}{\sigma}

\newcommand{\F}{\Phi}
\newcommand{\G}{\Gamma}

\renewcommand{\O}{\Omega}

\font\tenrm=cmr10

\font\bss=cmssdc10 at 12 pt 
\font\bssl=cmssdc10 at 16 pt
\font\cmsslll=cmss10 at 14 pt

\newcommand{\bC}{\mathbb{C}}
\newcommand{\bR}{\mathbb{R}}

\newcommand{\bN}{\mathbb{N}}
\newcommand{\bO}{\mathbb{O}}

\renewcommand{\gg}{\mathfrak{g}}

\newcommand{\gl}{\mathfrak{l}}

\newcommand{\gX}{\mathfrak{X}}
\newcommand{\gS}{\mathfrak{S}}

\newcommand{\gso}{\mathfrak{so}}

\newcommand{\gspin}{\mathfrak{spin}}
\newcommand{\gsp}{\mathfrak{sp}}

\newcommand{\ggl}{\mathfrak{gl}}

\newcommand\Sp[1]{\mathrm{Sp}(#1)}

\newcommand\SO[1]{\mathrm{SO}(#1)}
\newcommand\SU[1]{\mathrm{SU}(#1)}
\newcommand\U[1]{\mathrm{U}(#1)}
\newcommand\Spin[1]{\mathrm{Spin}(#1)}

\newcommand{\cb}{\mathcal{B}}

\newcommand{\cd}{\mathcal{D}}

\newcommand{\R}{\mathcal{R}}

\newcommand{\cu}{\mathcal{U}}

\newcommand\fr[2]{\tfrac{#1}{#2}}
\newcommand{\p}{\partial}

\renewcommand{\square}{\kern1pt\vbox
               {\hrule height 0.6pt\hbox{\vrule width 0.6pt\hskip 3pt
    \vbox{\vskip 6pt}\hskip 3pt\vrule width 0.6pt}\hrule height0.6pt}
    \kern1pt}

\newcommand{\La}{\wedge}
\newcommand{\ra}{\rightarrow}
\DeclareMathOperator\tr{tr\;}

\DeclareMathOperator\End{End\;}

\DeclareMathOperator\Hol{Hol\;}
\newcommand\der[1]{\tfrac{\partial}{\partial #1}}

\newcommand{\wt}{\widetilde}
\newcommand{\wh}{\widehat}
\newcommand{\ol}{\overline}

\newtheorem{Th}{Theorem}
\newtheorem{Prop}{Proposition}
\newtheorem{Cor}{Corollary}
\newtheorem{Lem}{Lemma}
\newtheorem{Def}{Definition}

\newcommand{\bt}{\begin{Th}\ \ }
\newcommand{\et}{\end{Th}}
\newcommand{\bp}{\begin{Prop}\ \ }
\newcommand{\ep}{\end{Prop}}
\newcommand{\bc}{\begin{Cor}\ \ }
\newcommand{\ec}{\end{Cor}}
\newcommand{\bl}{\begin{Lem}\ \ }
\newcommand{\el}{\end{Lem}}
\newcommand{\bd}{\begin{Def}\ \ }
\newcommand{\ed}{\end{Def}}

\newcommand{\pf}{\noindent{\it Proof:\ \ }}
\newcommand{\qed}{\hfill\square}
\newcommand{\n}{\nabla}
\newcommand{\op}{\oplus}
\newcommand{\ot}{\otimes}
\newcommand{\dpp}{\p_{++} }
\newcommand{\App}{A_{++} }
\newcommand{\npp}{\n^S_{\dpp} }
\newcommand{\dmm}{\p_{--} }
\newcommand{\Amm}{A_{--} }
\newcommand{\nmm}{\n^S_{\dmm} }
\newcommand{\dpm}{\p_{\pm\pm} }
\newcommand{\Apm}{A_{\pm\pm} }
\newcommand{\npm}{\n^S_{\dpm} }
\newcommand{\nz}{\n^S_{\p_0} }
\newcommand{\bop}{\bigoplus}

\newcommand{\Fone}{\stackrel{(1)}{F}{}^{(ee')}}
\newcommand{\Ftwo}{\stackrel{(2)}{F}{}^{[ee']}}
\newcommand{\Fthree}{\stackrel{(3)}{F}{}^{(ee')}}

\newcommand{\be}{\begin{equation}}
\newcommand{\ee}{\end{equation}}
\newcommand\la[1]{\label{#1}}
\newcommand\re[1]{\eqref{#1}}
\newcommand{\arr}{\begin{array}{rlll}}
\newcommand{\ea}{\end{array}}
\newcommand{\bea}{\begin{eqnarray}}
\newcommand{\eea}{\end{eqnarray}}
\newcommand{\bean}{\begin{eqnarray*}}
\newcommand{\eean}{\end{eqnarray*}}

\numberwithin{equation}{section}

\def\hk{hyper-K\"ahler}
\def\qk{quaternionic K\"ahler manifold}

\begin{document}
\begin{titlepage}
\rightline{math.DG/0209124}
\vskip 1.5 true cm
\begin{center}
{\bssl Yang-Mills connections over manifolds with Grassmann structure}

\vskip 1.0 true cm
{\cmsslll
Dmitri V.\ Alekseevsky$\,^{1}$, \ Vicente Cort\'es$\,^{2,3}$,\
Chandrashekar Devchand$\,^{2}$}
\vskip 0.8 true cm

{\tenrm {$^1$ Dept. of Mathematics, University of Hull,
Cottingham Road, Hull, HU6 7RX, UK}\\
{  $^2$ Mathematisches Institut der Universit\"at Bonn,
Beringstra\ss e 1, D-53115 Bonn, Germany}\\
{  $^3$  Institut \'Elie Cartan, 
Universit\'e Nancy 1, 
B.P. 239, F-54506 Vandoeuvre-l\`es-Nancy Cedex, France}\\
{D.V.Alekseevsky@maths.hull.ac.uk}\ ,\
{ cortes@iecn.u-nancy.fr}\ ,\
{  devchand@math.uni-bonn.de}}
\end{center}
\vskip 0.8 true cm

\begin{abstract}

\noindent
Let $M$ be a manifold with Grassmann structure, i.e.\  with an isomorphism of
the cotangent bundle $T^*M{\cong}E{\otimes}H$ with the tensor product of two
vector bundles $E$ and $H$. We define the notion of a half-flat connection
$\nabla^W$ in a vector bundle $W{\ra} M$ as a connection whose curvature 
$F\in S^2E\otimes\wedge^2 H\otimes\End W \subset\wedge^2 T^*M\otimes\End W$. 
Under appropriate assumptions, for example, when the Grassmann structure is
associated with a quaternionic K\"ahler structure on $M$, half-flatness
implies the Yang-Mills equations. Inspired by the harmonic space approach,
we develop a local construction of (holomorphic) half-flat connections
$\nabla^W$ over a complex manifold with (holomorphic) Grassmann structure
equipped with a suitable linear connection. Any such 
connection $\nabla^W$ can be
obtained from a prepotential by solving a system of linear first order ODEs. 
The construction can be applied, for instance, to the complexification of
hyper-K\"ahler manifolds or more generally to hyper-K\"ahler manifolds with
admissible torsion and to their higher-spin analogues. It yields solutions
of the Yang-Mills equations.

\end{abstract}
\vfill \hrule width 3.0 cm
{\small \noindent This work was supported by the `Schwerpunktprogramm
 Stringtheorie' of the Deutsche For\-schungs\-ge\-mein\-schaft, 
MPI f\"ur Mathematik (Bonn) and the SFB 256 (University of Bonn).}
\end{titlepage}

{\small 
\tableofcontents}
\vfil
%%%%%%%%%%%%%%%%%%%%%%%%%%%%%%%%%%%%%%%%%%%%%%%%%%%%%%%%
\section{Introduction}
The Yang-Mills self-duality equations have played an important role in field
theory and in differential geometry. They are the main source of examples of
solutions of the Yang-Mills equations on four-dimensional manifolds \cite{AHS}. 
The self-duality equations $\,*F^\n {=} F^\n\,$ mean that the curvature
$F^\n$ of a connection $\n$ over a Riemannian 4-fold $M$ is an eigenvector
of the Hodge star operator, associated with the volume 4-form, which acts on
two-forms. This apparently four-dimensional construction has an analogue in
Riemannian manifolds $M$ of arbitrary dimensions. Any 4-form $\Omega$ on $M$
defines an endomorphism $B_\Omega$ of the space of 2-forms and one can define
$(\O,\l)$-self-duality as the condition, $\,B_\O F^\n{=} \l F^\n$, that the
curvature is an eigenvector of $B_\O$ with eigenvalue $\,\l ={\rm const}\neq
0$. Under appropriate assumptions on $\Omega$ (for example, if it is
co-closed) this implies the Yang-Mills equations, just as in four
dimensions. For instance, this works for a constant $\Omega$ in flat space
\cite{CDFN,W} and for a parallel 4-form on a Riemannian manifold with special
holonomy (some examples are discussed in \cite{CS,N,BKS,DT,T}). If $\Omega$
is, for example, the canonical parallel 4-form associated to a quaternionic
K\"ahler manifold $M$ of dimension $4m$, then the eigenspaces of 
$B_\O$ are the irreducible
$\Sp{m}\cdot \Sp{1}$--submodules of the space of two-forms. In terms of the
associated locally defined Grassmann structure $T^{*\bC}M = E \otimes H$,
i.e.\  the identification of the complexified cotangent bundle $T^{*\bC}M$
with a tensor product of two vector bundles $E$ and $H$ of rank $2m$ and $2$
respectively, the eigenspace decomposition is given by 
\[
\wedge^2 T^{*\bC} M = 
S^2 E \ot \wedge^2 H\ \op\ \wedge_0^2 E\ot S^2 H\  \op\ \o_E\ot S^2 H\   ,
\] 
with corresponding
$B_\O$--eigenvalues $\l_1{=}1,\,\l_2{=}-1/3,\,\l_3{=}-(2m{+}1)/3$ \cite{W,CGK}. 
Here $\o_E$ and $ \o_H$ are two-forms on $E^*$ and $H^*$ such that the 
complex metric on $T^\bC M$ is given by $\o_E\ot \o_H$ and  $\wedge_0^2 E$
denotes the traceless part of $\wedge^2 E$ with respect to $\o_E$. 
The eigenspaces of $B_\O$ can thus be described in terms of the Grassmann 
structure, which is a natural generalisation of the
well-known spinor decomposition of a vector in four dimensions. 
A 2-form on any manifold with Grassmann structure is called
half-flat if it belongs to the eigenspace $S^2E\ot \wedge^2H$ and a
connection $\n$ with half-flat curvature is called half-flat. If the
Grassmann structure is associated with the quaternionic K\"ahler structure,
then a half-flat connection is the same as an  $(\O,\l_1)$-self-dual
connection and hence satisfies the Yang-Mills equations. Inspired by the
harmonic space approach \cite{GIOS}, we develop a construction of
locally-defined holomorphic half-flat connections on a manifold $M$ with
holomorphic admissible half-flat Grassmann structure, namely, a holomorphic
Grassmann structure $T^*M=E\ot H$ with holomorphic connections $\n^E$ and
$\n^H$ in the bundles $E$ and $H$ respectively,  such that $\n^H$ is flat and
the torsion of the linear connection  $\n=\n^E\ot {\rm Id}+ {\rm Id}\ot \n^H$
has no component in $S^3H\ot E^*\ot \wedge^2E$. The construction associates
to a holomorphic prepotential a half-flat connection through the solution of 
a system of linear first order ODEs. The construction can be applied, for
instance, to the complexification of hyper-K\"ahler manifolds or, more
generally, to hyper-K\"ahler manifolds with admissible torsion. Our
construction of gauge fields on such curved backgrounds extends that of
\cite{GIOS}, where flat torsionfree backgrounds were considered.
Moreover, we provide a geometric description of the harmonic space method 
of \cite{GIOS}.

We note that using analytic continuation any real analytic connection $\n$
over a real analytic Grassmann manifold allows extension to a holomorphic
connection $\n^\bC$ over a holomorphic Grassmann manifold and $\n$ can be
reconstructed from $\n^\bC$ in terms of some antiholomorphic involution. 

The main idea of our construction is to pull-back a half-flat connection
$\n$ in a holomorphic vector bundle $\nu:W\ra M$ to the harmonic space $S_H$. 
The latter is the space of all symplectic frames $h=(h_+,h_-)$ in the vector
bundle $H^*$. The group $\Sp{1,\bC}$ acts freely on $S_H$, with the orbit 
space $S_H/\Sp{1,\bC}=M$. Hence, the projection $\pi:S_H\ra M$ is an
$\Sp{1,\bC}$-principal bundle. Choosing a (local) trivialisation,
$M \ni x \mapsto (h_1(x),h_2(x)) \in S_H$, of $\pi$ we can make the
identification $S_H=\Sp{1,\bC}\times M$. There exists a canonical
decomposition,   
$$ TS_H = T^vS_H \op \cd_+ \op \cd_-\ , $$ 
of the (holomorphic) tangent bundle into the 
vertical sub-bundle $T^vS_H $ and two
(holomorphic) distributions $\cd_+$ and $\cd_-$ spanned respectively by vector
fields $X^e_+$ and $X^e_-$ canonically associated with sections $e$ of the 
bundle $E^*$. If the Grassmann structure is admissible and half-flat,
the distributions $\cd_+$ and $\cd_-$ are integrable. The vertical 
distribution $T^vS_H$ is spanned by vector fields $\p_0, \dpp, \dmm$, which
correspond to the standard generators of the Lie algebra $\gsp(1,\bC)$. 
A half-flat connection $\n$ in the bundle $\nu:W\ra M$ induces the pull-back
connection $\pi^* \n$ in the pull-back bundle $\pi^* \nu: \pi^*W \ra S_H$. 
Since $\n$ is half-flat, the curvature $F$ of $\pi^* \n$ satisfies certain 
equations (see Definitions \ref{half-flatDef} and \ref{half-flatSDef}).  
A connection in $\pi^* \nu$ satisfying these equations is called a half-flat
connection over $S_H$ and is gauge equivalent to the pull-back of a half-flat
connection over $M$. Any half-flat connection over $S_H$ is flat along the
leaves of the integrable distribution $\langle \cd_+, \p_0\rangle $ spanned
by $\p_0$ and $\cd_+$. We can therefore choose a frame of the vector bundle
$\pi^* \nu$ which is parallel along its leaves. Such a frame is called an
analytic frame. With respect to such a frame a half-flat connection has no
potentials in the directions of the distribution  $\langle\cd_+,\p_0\rangle
$.  Starting from a matrix-valued function (prepotential) $\App$ on $S_H$,
which is constant along the leaves of the distribution $\cd_+$ and satisfies
the homogeneity condition $\p_0 \App = 2 \App\,$, we construct a connection
which satisfies almost all the conditions of half-flatness. We call such a
connection an almost half-flat connection. It is half-flat if and only if its
curvature satisfies the equation $F(\dmm,\cd_-)=0$. The construction of an
almost half-flat connection reduces to the solution of first order linear
ODEs. Assuming that the almost half-flat connection $\n$ is defined globally
along the fibres (over $\pi^{-1}U$, where $U\subset M$ is a domain in $M$) we
can modify $\n$ to a half-flat connection over $S_H$ which is the pull-back
of a half-flat connection over $M$. In order to do this, we  re-write $\n$
with respect to a ``central frame'', namely,  a frame  parallel along the
fibres of $\pi$. The transformation from the analytic to  the central frame
reduces to the solution of the system of equations  $$ 
\dpp \F = - \App \F\quad,\quad \p_0 \F =0\ . 
$$
With respect to the central frame the potential $C(X^e_+)$ of the connection
$\n$ in the direction of the vector field $X^e_+\in \cd_+$ has the form 
$C(X^e_+)= u_+^\a C^e_\a$, where $C^e_\a$ are matrix-valued functions on
$M= M \times \{{\rm Id}\} \subset M \times \Sp{1,\bC}$ and $u_\pm^\a\ ,\
\a=1,2$, are matrix coefficients of $\Sp{1,\bC}$. The matrix-valued functions 
$C^e_1, C^e_2$ define the desired half-flat connection on $M$ given by
$$
\n^M_{e\ot h_1} = e\ot h_1 + C^e_1 \quad,\quad
\n^M_{e\ot h_2} = e\ot h_2 + C^e_2\ . 
$$
Moreover, any half-flat connection may be obtained in this way. 

The above construction allows generalisation to  manifolds with spin
$\fr{m}{2}$ Grassmann structure. This means that the cotangent bundle is
identified as $T^{*}M = E\ot F = E \otimes S^mH$, where $E$ and $H$ are
(holomorphic) vector bundles of rank $p$ and 2 respectively. If a connection
$\n^E$ on $E$ and a flat connection $\n^H$ on $H$ are given, then the
Grassmann structure is called half-flat. The connection $\n^H$ defines a flat
connection $\n^F$ on $F=S^mH$ and the linear connection 
$\n=\n^E\ot{\rm Id} + {\rm Id}\ot \n^F$. The associated harmonic
space $\pi:S_H\ra M$ is defined as above, as the space of all symplectic
frames $h=(h_+,h_-)$ in $H^*$. Its tangent space has decomposition
$$ 
TS_H = T^vS_H \op \bigoplus_{k=0}^{m} \cd_{k+} \op 
\bigoplus_{k=1}^{m} \cd_{k-}\ . 
$$
Under certain conditions on the torsion of $\n$ the distribution
$\cd^k_{(+)} := \bop_{i=0}^{k}\cd_{(m-2i)+}\,$, $k\le m/2$, is integrable. 
Such a half-flat Grassmann structure is called $k$-admissible. Generalising 
the notion of a half-flat connection, we may define a $k$-partially flat
connection $\n$ over a manifold with half-flat spin~$\fr{m}{2}$ Grassmann
structure such that the pull-back connection $\pi^* \n$ has no curvature in
the directions of $\cd^k_{(+)}$. The harmonic space method can be applied to
construct $k$-partially flat connections over $k$-admissible half-flat 
spin~$\fr{m}{2}$ Grassmann manifolds. In the final section we consider the
case of $m=3$ and sketch the construction of 0-- and 1--partially flat 
connections. The latter are Yang-Mills connections. 

We should like to thank Andrea Spiro for useful comments on the manuscript.

\section{Generalised self-duality for  manifolds of dimension greater
than four}

\subsection{Yang-Mills data}

Let  $\nu : W \rightarrow M$  be a real vector bundle over $M$ 
and $\n$ a connection in $\nu$, that is a bilinear map 
$$\arr
\n : \gX(M) \times \G(\nu) &\ra& \G(\nu)\\
           (X,  \s )&\mapsto& \n_X \s\ 
\ea$$
which is $C^\infty(M)$-linear in the vector field $X\in\gX(M)$ and
satisfies the Leibniz rule $ \n_X(f\s) = (Xf) \s + f \n_X \s $, 
for any function $f\in C^\infty(M)$ and any section,  $\s\in \G(\nu)$,
of $\nu$. The map $\n$ can be extended to a complex bilinear map,
\be\arr
\n : \gX^\bC(M) \times \G(W^\bC{\ra}M) &\ra& \G(W^\bC{\ra}M)\\
           (X,\s ) &\mapsto& \n_X \s\ ,
\ea\la{cnab}\ee
where $\gX^\bC(M)$ is the space of complex vector fields
$X{+}iY\,;\, X,Y{\in}\gX(M)$ and $ W^\bC \rightarrow M$ is the
complexification of the vector bundle $\nu$. Note that $\n$
satisfies the reality condition
\be
\n_{\ol X} \ol\s = \ol{\n_X \s}\ ,\
X\in \gX^\bC(M)\,,\,\s\in\G(W^\bC{\ra}M) \, ,
\la{reality}
\ee
where the bar denotes complex conjugation. Conversely, any
$\bC$-bilinear map \re{cnab} which is $\gX^\bC(M)$-linear and
satisfies the Leibniz rule and the reality condition \re{reality},
defines a connection $\n$ in the real vector bundle $\nu$. If the
reality condition \re{reality} is dropped, then \re{cnab} defines a
connection in the complex vector bundle $W^\bC{\ra}M$.

Let $\f=(\f_1,\dots,\f_r)$ denote a local frame of $\nu$ such that
for any section $\s\in \G(\nu)$, $\s=\sum s^i \f_i =: \f \cdot s $,
where $s^i$ are the coordinates of $\s$ with respect to the frame $\f$ 
and $s = (s^1, \dots,  s^r)^t$. 
Then the connection $\n$  in $\nu$ has local expression
$$
\n_X \s = \n_X (s^i \f_i) =  \f \cdot \n_X s :=
(X s^i + \sum_j A^i_j(X) s^j ) \f_i\, ,
$$
where $ A^i_j(X) = (\n_X \f_j , \f^i) $ and $\f^*=(\f^1,\dots,\f^r)$ 
denotes the dual frame. The locally defined matrix-valued
1-form $A=(A^i_j): M \ra \ggl(r,\bR)$
is called the Yang-Mills potential with respect to the frame $\f$. 
If the vector bundle $\nu$ has structure group $G$, i.e.\  if it is a
bundle associated with a principal $G$-bundle  $P\ra M$ and
a representation $\rho: G\ra GL(r,\bR)$, such that $W= P\times_G \bR^r$,
then we may always choose a frame $\f$ for which the potential takes values
in the Lie algebra $\gg = {\rm Lie}\ \rho(G) \subset \ggl(r,\bR)$. 
We will symbolically write  $\n_X = X + A\ ,\ A = A^\f $. 
A change of frame (gauge transformation) $\f'= \f U$ induces changes
$s' =U^{-1} s $ and
$\f' (X+A'(X))s'\ =\ \f'\n_X s'\ =\ \f\n_X s\ =\ \f(X+A(X))s\ =\
\f'U^{-1}(X+A(X))Us'$,
yielding the transformation rule for the potential,
\be
A' = U^{-1} (XU) +  U^{-1} A(X)U  = U^{-1} \n_X U\, . 
\la{gt}\ee
The curvature of the connection $\n$,
$ F = F^\n \in \O^2_M(\End W) = \G(\La^2 T^* M\ot \End W)$,
is given by
\bean
F(X,Y) &=& \left[ \n_X , \n_Y \right] - \n_{[X,Y]} \\ 
         &=&  XA(Y) - YA(X) + \left[ A(X) , A(Y) \right] - A([X,Y])\ . 
\eean
The Jacobi identity for $\n_X$ is equivalent to the Bianchi identity,
$\ d^\n F^\n = 0$. 
Here the covariant derivative $ d^\n : \O^p(\End W) \ra  \O^{p+1}(\End W)$
is defined by
$$
d^\n ( \o\ot C ) = d\o \ot  C\  +\  (-1)^p\o \wedge \n C\ ,
$$
where $\o$ is a $p$-form and $C$ is a section of $\End W$. (The connection 
$\n$ on $W$ induces a connection on $\End W $ denoted by the same symbol.)

On any $n$-dimensional oriented pseudo-Riemannian (or complex Riemannian)
manifold, $(M ,g)$ using the canonical volume form 
$\mbox{vol}^g \in \La^n T^*M$, we define the Hodge $\ast$ operator which
interchanges forms of complementary degree, $*: \La^p T^*M\ra \La^{n-p}T^*M$, 
by the relation: $\langle\a,\b\rangle\mbox{vol}^g=\a \wedge *\b\,$,
where $\a,\b \in \La^p T^*M$ and $\langle .\,,\,.\rangle $ is the natural
scalar product on $\La^p T^*M$ induced by the metric $g$. 
We define $\,\ast : \wedge^pT^*M{\ot}\End  W \ra \wedge^{n-p}T^*M{\ot}\End W$ 
by $\,*( \o\ot C ):=( *\o \ot C )$. 
\bd
Let $\nu: W\ra M$ be a real vector bundle over a pseudo-Riemannian manifold
$(M,g)$. A {\bss YM connection} $\n$  in $\nu$ is one which satisfies the
Yang-Mills equation
$$
d^\n \ast F^\n = 0\ . 
$$
\ed
On a closed manifold this is the Euler-Lagrange equation for the 
YM functional 
\be
|| F^\n ||^2 = \int_M\  \left|F^\n \right|^2\  {\rm vol}^g ,
\la{yma}\ee
where the norm on $\La^2 T^{\ast}M \otimes \End W$ is induced by the
pseudo-Riemannian metric on $M$ and the natural metric on $\End W$.

\subsection{Self-duality conditions}
On a Riemannian four-manifold, the $\ast $ operator maps 2-forms to 2-forms
and has eigenvalues $\pm1$. The curvature tensor therefore has decomposition
into the eigenspaces of the $\ast $ operator,
$$
F^\n
     = F^\n_{+1} \op F^\n_{-1} \in \O^+_M(\End W) \op \O^-_M(\End W). 
$$
This splitting corresponds to the decomposition of
the $\SO4$-module $\La^2 \bR^4 =
\La^2_+ \op  \La^2_- \cong \gso(4)=\gsp(1)\op\gsp(1)$ into
its irreducible submodules. 
We call $\n$ and $F^\n $  self-dual or anti-self-dual if
$F^\n_{-1}:= \fr12 (F^\n-\ast F^\n)=0$ or
$F^\n_{+1}:= \fr12 (F^\n+\ast F^\n)=0$, respectively. 
For (anti-) self-dual connections,
the YM equation, $\,d^\n \ast F^\n = 0$, is an immediate consequence of the
Bianchi identity, $\,d^\n F^\n = 0$. On  closed manifolds (anti-) self-dual
connections in fact minimise the YM functional \re{yma}, since the inequality
$$
     || F^\n ||^2 = || F^\n_{+1} ||^2 + || F^\n_{-1} ||^2 \ge
\left|{|| F^\n_{+1} ||^2 - || F^\n_{-1} ||^2}\right|
= 8\pi^2 \left|c_2(W)[M]\right|
$$
is saturated. Here $c_2(W)[M] = \fr{1}{8\pi^2}\int_M \tr \ F^\n\wedge F^\n$
is the evaluation of the second Chern class of the bundle $W$ on the
fundamental cycle. 

The apparently four-dimensional notion of self-duality, has an analogue
in higher dimensions. 
The construction originally given in \cite{CDFN}
for flat spaces, extends to arbitrary manifolds $(M ,g)$, of dimension
greater than four, as follows. 

For $\O \in \O^4(M)$ we define a symmetric tracefree 
endomorphism field $B_\O : \La^2 T^*M \ra \La^2 T^*M$ by
\be
B_\O \o := \ast  (\ast  \O \wedge \o )\ ,
\la{Bdef}\ee
where $\o\in \La^2 T^*M$. 
This endomorphism is zero if and only if the 4-form $\O$ is zero. 
Moreover, we have:
\bl \label{contractionLemma} Let
\[ \O = \sum \O_{ijkl}e^i \wedge e^j \wedge e^k \wedge e^l, \quad \o = 
\sum \o_{ij}e^i \wedge e^j\]
be the expressions for $\O$ and $\o$ with respect to a frame $e^i$ of $T^*M$. 
Then $B_\O$ is given as the contraction
\[ B_\O\o 
= 12 \sum g^{ii'}g^{jj'}\ \O_{ijkl}\ \o_{i'j'}\ e^k \wedge e^l \, . 
\]
\el

\pf It is sufficient to check the above formula
for decomposable forms $\O = e^i\wedge e^j\wedge e^k
\wedge e^l$ and $\o = e^m\wedge e^n$, where the $e^i$ form an
orthonormal basis of $T^*M$. \qed

\bd
A 4-form $\ \O \in \O^4(M)\,$ on a pseudo-Riemannian manifold $M$ is
called {\bss appropriate} if there exists a nonzero real constant
eigenvalue $\l$ of the endomorphism field  $B_\O$. 
\la{approp}
\ed
We note that on a Riemannian manifold the eigenvalues of $B_\O$ are
real for any 4-form $\O$. 
A generalisation of the four dimensional notion of self-duality
may now be defined:

\bd\la{sdy}
Let $\O$ be an appropriate 4-form on a pseudo-Riemannian manifold $(M,g)$
and $\l\neq 0 \in \bR$. A connection $\n$ in a vector bundle $\nu: W\ra M$
is ${\boldsymbol{(\O , \l )}}${\bss{-self-dual}} if its curvature $F^\n$
satisfies the linear algebraic system
\bea
B_\O F^\n &=& \l F^\n \la{sdy1}\\
(d\ast \O) \wedge F^\n &=& 0\ . 
\la{sdy2}
\eea
\ed

\bt \label{YMThm}
Let $(M,g)$ be a pseudo-Riemannian manifold with an appropriate 4-form $\O$. 
Then any $(\O , \l )$-self-dual connection  $\n$ is a YM connection. 
\et

\pf  
Using \re{sdy1} and \re{Bdef} we obtain
\bean
d^\n * F^\n &=&  \fr{1}{\l} d^\n * B_\O F^\n \
=\  \pm\fr{1}{\l} d^\n (*\O\wedge F^\n)  \\ 
&=& \pm\fr{1}{\l}\left((d *\O)\wedge F^\n + *\O\wedge d^\n F^\n \right)\\
&=& 0 \, ,
\eean
in virtue of \re{sdy2} and the Bianchi identity $d^\n F^\n{=}0$.
\qed

\noindent
Examples of manifolds admitting appropriate 4-forms are easily obtained. 
Let $V$ be a pseudo-Euclidean vector space and
$G \subset \SO{V}$ be a linear group preserving a nonzero
element $\O_0 \in \La^4 V$. Denote by $\O_{ijkl}$ the components of
$\O_0$ with respect to an orthonormal basis of $V$. Given a manifold
$M$ with a $G$-structure, $\pi: P\ra M$, i.e.\  a principal $G$-subbundle of 
the bundle of frames on $M$, we can define a 4-form
$\O := \sum \O_{ijkl}e^i\wedge e^j\wedge e^k\wedge e^l$,
where $(e^1, \ldots , e^n)$ is a coframe dual to a $G$-frame
$p=(e_1, \ldots , e_n)\in P$. Since $G \subset \SO{V}$,
$M$ has the structure of an oriented pseudo-Riemannian manifold and
we can define the operator $B_\O$. The matrix components
of $B_\O = \sum B_{ij}^{kl}e^i\wedge e^j\otimes e_k\wedge e_l$
are constant for any $G$-frame and so are its eigenvalues. Hence $\O$ is
appropriate if the endomorphism $B_{\O_0}\in \La^4 V$ has a nonzero
real eigenvalue $\l$. This is automatic in the Riemannian case.

There exist many examples of subgroups $G\subset\SO{V}$ admitting nonzero
$G$-invariant 4-forms, as shown by the following construction. Let 
$G\subset \SO{V}$ be a closed subgroup of the pseudo-orthogonal group
$\SO{V}$ and $\gg{\subset} \gso(V){\cong} \La^2 V^*$ its Lie algebra. Assume
that $\gg$ admits a $G$-invariant symmetric nondegenerate bilinear form
$B\in S^2(\gg^*)^G$, where $W^G$ denotes the space of $G$-invariant
elements of a $G$-module $W$. 
We can then identify $\gg$ with its dual $\gg^*$ via $B$ and consider
$B$ as an element of $(S^2(\gg))^G {\subset} (S^2 \La^2 V^*)^G$. 
A $G$-invariant 4-form is then defined by
$\,\O_0^G {:=} \mbox{alt\,} B \in (\La^4 V^* )^G$,
where $\mbox{alt\,}{:} S^2\La^2V^* {\rightarrow} \La^4 V^*$ denotes 
alternation. 
We denote the corresponding 4-form on a manifold with $G$-structure by $\O^G$. 
The following variant of a theorem by Kostant \cite{K} provides a wealth of
examples of nonzero $\O_0^G\,$'s. 

\bt
Let $G \subset \SO{V}$ be a closed subgroup whose Lie algebra $\gg$
admits a nondegenerate $G$-invariant bilinear form $B\in (S^2\gg)^G$. 
If the $G$-module $V$ is not equivalent to the isotropy module of a
pseudo-Riemannian symmetric space, then the four-form
$\O_0^G := \mbox{alt\,} B \in (\La^4 V)^G$ is nonzero. 
\label{KostantThm}
\et

\pf
Recall that the $\SO{V}$-module $ S^2\La^2 V$ decomposes according to
$\  S^2\La^2 V = \R(\gso(V)) + \La^4 V $, where
$\R(\gso(V))$ denotes the space of curvature tensors of type  $\gso(V)$,
i.e.\  the space of two-forms fulfilling the first Bianchi identity or
the kernel of the map $\;\mbox{alt\,}: S^2\La^2V \rightarrow \La^4 V$. 
If $\O_0^G = \mbox{alt\,} B =0$, then $B$ is a nonzero element of
$ \R(\gso(V))\cap S^2(\gg)^G = \R(\gg)^G$. 
Since $B$ is a $G$-invariant two-form on $V$ with values in $\gg$ it
can be used to define a Lie bracket $[\,\cdot\,,\,\cdot\,]$ on the vector
space  $\gl=\gg \op V$ thus:\\
(i)~ $\gg$ is a subalgebra of $\gl$,\\
(ii) $V$ is a $\gg$-submodule with action defined by the
inclusion $\gg \subset \gso(V)$,\\
(iii) $[u,v] := B(u,v) \in \gg$ if $u,v \in V$. \\
The Jacobi identity follows from the Bianchi identity and the $G$-invariance. 
Let $L$ be the simply connected Lie group with Lie algebra $\gl$. Then
$L/G_0$ is a Riemannian symmetric space with  $V$ as its isotropy module,
where $G_0 \subset L$ is the connected Lie subgroup with
${\rm Lie\,} G_0 = \gg$. 
\qed

Clearly, \re{sdy2} is automatic if the 4-form $\O$
is  co-closed, $d\ast \O=0$. This is the case, for example, if $\O$ is
parallel. In the Riemannian case the Berger list of irreducible holonomy
groups \cite{B} and a theorem of Kostant \cite{K} yield the
following result. 

\bt
Let $M$ be a complete simply connected irreducible
Riemannian manifold of dimension $n\ge 4$ with holonomy group
$\Hol\subset \SO{n}$, $\Hol\neq \SO{n}$. 
Then  $M$ admits a nontrivial parallel
4-form if one of the following holds:
i) $M$ is not a symmetric space or\\
ii) $M$ is a symmetric space and has a non-simple holonomy or,
equivalently, isotropy group. 
\et

\pf By Berger's theorem on Riemannian irreducible holonomy groups \cite{B}, 
we have:\\
a) $M$ is not a symmetric space and its holonomy group is one of:
$\U{\fr{n}{2}}$, $\SU{\fr{n}{2}}$, $\Sp{\fr{n}{4}} \Sp1$, $\Sp{\fr{n}{4}}$,
${\rm G_2}$, $\Spin7$,
or\\
b) $M$ is a symmetric space. \\
All the groups in a) admit invariant 4-forms. These are given below. 
A theorem of Kostant \cite{K} states that a simply-connected
irreducible Riemannian symmetric space $G/K$ has no nonzero parallel
4-form if and only if  the isotropy group $K$ is simple. 
\qed

In the following examples we explicitly describe parallel (hence appropriate)
4-forms $\O$ on Riemannian $n$-manifolds with  holonomy groups
$\Hol \neq \SO{n}$ from Berger's list. 

\noindent{\bf 1.\ } K\"ahler manifolds,
$\Hol \subset {\mathrm U(m)} \subset \SO{2m}$,
$n=2m$:   $\O=  \o \wedge \o $, where $\o$ is the K\"ahler form. 
One can check that this is proportional to $\O^{\SU{m}}$
and that any parallel 4-form is proportional to $ \o \wedge \o $  if
the holonomy group is $\SU{m}$ or ${\mathrm U(m)}$. 
If $\Hol\subset\Sp{k} \subset \SU{2k}\subset \SO{4k}$, $n=4k > 4$,
i.e.\  if the manifold is \hk, there exist three skewsymmetric parallel 
complex structures $J_\a, \a=1,2,3$. Then there exist six independent
parallel 4-forms $\o_\a \wedge\o_\b,\ \a,\b=1,2,3$, where $\o_\a$ is the
K\"ahler form associated to $J_\a$. For low dimensional examples,
eigenvalues and eigenspaces of $B_\O$  are given in \cite{CDFN}. 

\noindent{\bf 2.\ } Quaternionic K\"ahler manifolds,
$\Hol\subset\Sp{m} \Sp1 \subset \SO{4m}$, $n=4m$. In this case
there exist 3 locally defined almost complex structures  $J_\a$,
with corresponding K\"ahler forms $\o_\a$, such that
the 4-form $\O:=\sum_\a \o_\a \wedge \o_\a$ is globally defined and parallel. 
This will be discussed in more detail in section \ref{qksubsection}. 

\noindent{\bf 3.\ } $\Hol\subset G_2 \subset \SO{7}$. Let $V = \bO =
\bR 1 + {\rm Im} \bO =  \bR \oplus \bR^7 = \bR^8$ be the algebra of
octonions. Recall that $G_2 = {\rm Aut}(\bO )$ is the group of
automorphisms of the octonions. We can decompose the product of two
octonions $a,b$ into its real and imaginary parts as follows:
\[ ab = a\cdot b = \langle a, b\rangle 1 + \fr{1}{2}[a,b] \, ,\]
where $\langle a, b\rangle$ is the scalar product  and $[a,b] = ab-ba$
is the commutator. We define a 3-form $\varphi$ and a 4-form
$\psi$ on ${\rm Im} \, \bO = \bR^7$ by the formulas
$$\arr
\varphi (x,y,z) &:=& \langle x\cdot y, z\rangle = \fr12
\langle [x,y], z\rangle \\[2pt] 
\psi (x,y,z,w) &:=& \langle [x,y,z], w \rangle \, ,
\ea$$
where $[x,y,z] = (xy)z - x(yz)$ is the associator. It is known that
$\psi = \ast \varphi$. Notice that $G_2$ is the group of isometries
of $\bO = \bR^8$ which fix the identity element $1$ and preserve the
3-form $\varphi$ (or equivalently the 4-form $\psi$) on ${\rm Im}\, \bO$. 
The 4-form $\psi$ defines a parallel 4-form on any Riemannian 7-fold
with holonomy $G_2 \subset \SO{7}$. It is known \cite{OV} that 
$\wedge^4 \bR^7 = \bR \psi \oplus V^7(\pi_1) \oplus V^{27}(2\pi_1)$, 
where $V^d(\pi )$ is the $d$-dimensional real irreducible representation of
$G_2$ with highest weight $\pi$ and $\pi_i$ denotes the $i$-th fundamental
weight of $G_2$. From this it follows that the 4-form
$\O_0^{G_2}$ coincides with $\psi$ up to scaling. The corresponding
endomorphism $B_\psi$ of $\wedge^2 \bR^7 = \gg_2 \oplus \bR^7$ has
two distinct eigenvalues which correspond to the two irreducible
$G_2$-submodules $\gg_2 $ and  $\bR^7 \subset \wedge^2 \bR^7$ \cite{CDFN}. 

\noindent{\bf 4.\ } $\Hol\subset \Spin7  \subset \SO{8}$. 
Using the 3- and 4-forms $\varphi$ and $\psi$ on $\bR^7$ introduced in
the $G_2$-case, we construct the 4-form
\[ \O = dt \wedge \varphi + \psi \, ,\]
where $t$ is the first coordinate on $\bR^8 = \bR 1 + \bR^7$. In particular,
$$
 \O (1,x,y,z) = \varphi (x,y,z) \quad,\quad
         \O (x,y,z,w) = \psi (x,y,z,w) \, ,\quad x,y,z,w \in \bR^7 \, . 
$$
This 4-form $\O$ defines a parallel 4-form on any Riemannian 8-fold
with holonomy $\Spin7  \subset \SO{8}$. It is known that $\wedge^4 \bR^8
= \bR \O \oplus V^7(\pi_1) \oplus V^{27}(2\pi_1) \oplus \wedge^4 \bR^7$. 
{}From this it follows that the 4-form $\O_0^{\Spin7}$ coincides with $\O$
up to scaling.  The corresponding endomorphism $B_\O$ of 
$\wedge^2 \bR^8 = \gspin_7 \oplus \bR^7$ has two distinct eigenvalues which
correspond to the two irreducible $\Spin7$-submodules $\gspin_7$ and  
$\bR^7\subset \wedge^2 \bR^8$ \cite{CDFN}.

\subsection{Quaternionic K\"ahler case}\la{qksubsection}
Now we discuss in more detail the case of \qk s (example 2 above). 
Riemannian manifolds $(M,g)$ with holonomy group $\Hol\subset\Sp{m}\Sp1$
are called \qk s. A \qk\ with holonomy group $\Hol\subset\Sp{m}$ is called
hyper-K\"ahler. 
On any \qk\ $M$, there exists a rank 3 vector subbundle $Q\subset \End TM$,
invariant under parallel transport, which
is locally spanned by three almost complex structures
$(J_\a) = (J_1, J_2, J_3{=}J_1 J_2 {=} - J_2 J_1)$. The latter are in general
only locally defined. The (globally defined)
vector bundle $Q$ is called the {\bss quaternionic structure} of $M$. 
A local frame $(J_\a)$ as above is called a {\bss standard frame} for $Q$. 
Similarly, a {\bss standard basis} of $Q$ at $m\in M$, is a triple
$I, J , K = IJ = - JI \in Q_m$ of complex structures on $T_mM$. 
A \qk\ is hyper-K\"ahler if and only if there exists a globally defined
parallel standard frame $(J_\a) = (J_1, J_2, J_3{=}J_1 J_2{=}- J_2 J_1)$. 

Given a standard frame, we may locally define three nondegenerate
2-forms $\o_\a := g(J_\a \cdot ,\cdot)$. The 4-form 
$$
\O:=\sum_\a \o_\a \wedge \o_\a \la{Omega}
$$
is independent of the choice of standard frame and defines a global parallel
4-form.

To describe the eigenspace decomposition of $\O$ it is convenient to use
the Grassmann structure (i.e.\  generalised spinor decomposition) of
a \qk . Recall that a {\bss Grassmann structure} on a (real) manifold $M$ is
defined as an
isomorphism $T^{*\bC} M \cong E\ot H$ of the complexified cotangent bundle
with the tensor product of two complex vector bundles $E$ and $H$ over $M$. 
Any \qk\ admits a (locally defined) Grassmann
structure $T^{*\bC}M = E \ot H$, where $H$ has rank 2, such that the
holonomy group ${\rm Hol} \subset \Sp{E} \ot \Sp{H}$. This
follows from the fact that any complex irreducible representation
of the group $\Sp{m} \times \Sp{1}$ is a tensor product of
irreducible representations of its factors. 

The complex extension $g^{\bC}$ of the Riemannian metric
defines a complex bilinear metric on $T^\bC M$,  which locally
factorises as $g^{\bC} = \o_E \ot \o_H$, where
$\o_E$ and $\o_H$ are sections of $\La^2E$ and
$\La^2 H$, defining complex symplectic forms on the fibres of $E^*$ and
$H^*$, respectively. We call $\o_E$ and $\o_H$
the symplectic forms of the symplectic vector bundles $E^*$ and $H^*$. 

In terms of the Grassmann structure the eigenspaces $V_{\l}$
of the endomorphism  $B_{\O}$ on $\La^2 T^{*\bC}M$ are
given by \cite{W,CGK}:
$$ V_{\l_1} = S^2 E\ot \o_H  \quad,\quad
                     V_{\l_2} = \La^2_0E\ot S^2H \quad,\quad
                     V_{\l_3} = \o_E\ot S^2H \, ,
$$
where $\La^2_0E$ is the space of $\o_E$-traceless 2-forms and
the eigenvalues are $\l_1 = 1$,  $\l_2 =-1/3$ and $\l_3 = -(2m+1)/3$. 
In particular  the $\l_1$-self-duality condition takes the form
\be F^{\n} \in  S^2 E\ot \o_H \ot \End W\ .
     \la{sdqk}
\ee
Note that since $\O$ is parallel it is appropriate and co-closed
and hence the $(\O ,\l)$-self-duality equations (Definition \ref{sdy})
reduce to \re{sdy1}, which implies the Yang-Mills equation. It is known 
(see Theorem 1 of \cite{CS}) that
$\l_1$- and $\l_3$-self-dual connections correspond to absolute
minima of the Yang-Mills functional on compact \qk s.

\subsection{Self-duality as half-flatness}
The $\l_1$-self-duality equation \re{sdqk} in fact depends only the
existence of the factorisation
$T^{*\bC}M \cong E\ot H$ and the symplectic structure in $H^*$. 
A connection $\n$ in a vector
bundle $W$ over a manifold $M$ with a Grassmann structure
is called {\bss half-flat} if its curvature satisfies the condition
\be
F^\n \in  S^2 E \ot \La^2 H \ot \End W \ . 
\la{sdgs}\ee
In general such half-flat connections are {\it not} YM connections
(with respect to some metric), but it is possible to impose further conditions
on $F^\n$ in order to enforce the YM equation. 
In fact it is the half-flatness of the connection,
rather than the YM property, which is crucial for our construction
of solutions. 
\bp \la{l1sdProp}
A connection $\n$ in a vector bundle $W \rightarrow M$ over a
quaternionic K\"ahler manifold is half-flat if and only if
it is $\l_1$-self-dual. Hence any such connection is a
Yang-Mills connection. 
\ep
\pf
The result follows from \re{sdqk} and \re{sdgs} since $\La^2H$ is the line
bundle generated by $\o_H$. 
\qed

\noindent The Levi-Civita connection on a hyper-K\"ahler manifold is an example
of a half-flat linear connection. Its complexification gives an example
of what we call an admissible half-flat Grassmann structure in the 
next section.  

\section{Manifolds with half-flat holomorphic Grassmann structure}

Our goal is to give a construction of half-flat connections in a
vector bundle $\nu:W\ra M$ over a manifold $M$. If all objects are
real analytic, using analytic continuation we may obtain
corresponding complex analytic objects. Specifically, assume that
the manifold $M$ and the bundle $\nu$ are real analytic. Then $M$
is defined by an atlas of charts with analytic transition
functions. Extending these functions to complex holomorphic
functions, we may extend $M$ to a complex manifold $M^\bC$ with
antiholomorphic involution $\tau$ such that $M= (M^\bC)^\tau$,
the fixed point set of $\tau$. Similarly, a real analytic vector
bundle $\nu: W\ra M$ can be extended to a holomorphic vector
bundle $\nu^\bC: W^\bC\ra M^\bC$. Moreover, an analytic connection
$\n$ in $\nu$ can be extended to a holomorphic connection $\n$ in
$\nu^\bC$. A holomorphic extension of a Yang-Mills connection is also a
Yang-Mills connection. In the rest of this paper, we shall assume that
all objects (manifolds, bundles and connections) are holomorphic. In section
\ref{construction} we shall give a construction of half-flat connections in 
a holomorphic bundle $W \rightarrow M$ over a complex manifold $M$ with
holomorphic Grassmann structure. Now we describe the required geometrical
notions. In particular, we provide a description of the harmonic spaces of 
\cite{GIOS} in geometric language. Our description affords application
to the construction of half-flat connections over more general manifolds than
the flat torsionfree backgrounds previously considered in the harmonic space 
literature (see e.g. \cite{GIOS}).

\subsection{Grassmann structure}

Let $M$ be a complex manifold with holomorphic Grassmann structure
$T^*M = E\ot H$, the isomorphism of the holomorphic cotangent bundle over 
$M$ with the tensor product of holomorphic vector bundles $E$ and $H$ over
$M$ of rank $p$ and $q$ respectively. Then $TM = E^*\ot H^*$.  
A holomorphic linear connection
$\n$ on $M$ is called a {\bss holomorphic Grassmann connection} 
if it preserves the holomorphic Grassmann structure. This means that for 
any vector field $X$ on $M$ and local sections $e\in\G(E)$ and $h\in 
\G(H)$,
$$
\n_X (e\ot h) = \n^E_X e \ot h + e \ot \n^H_X h
$$
where $\n^E,\n^H$ are connections in the bundles $E,H$ respectively. 
\bd
A holomorphic Grassmann structure, $T^*M = E\ot H$, on a complex manifold $M$
with a holomorphic Grassmann connection $\n=\n^E\ot{\rm Id}+{\rm Id}\ot \n^H$
is called {\bss half-flat} if the connection $\n^H$ in the holomorphic
vector bundle $H\ra M$ is flat. A manifold with such a half-flat holomorphic
Grassmann structure is called a {\bss half-flat Grassmann manifold}.  
\ed
{\bf Assumption. } In this section we assume that
$M$ is a manifold with a half-flat holomorphic Grassmann 
structure ($T^*M = E\ot H\ ,\n = \n^E\ot {\rm Id}+ {\rm Id}\ot \n^H$), 
such that $H$ has rank $2$ and that a $\n^H$-parallel nondegenerate fibre-wise
2-form $\o_H \in \G(\La^2 H)$ in the bundle $H^*$ is fixed. 
If, in addition, a $\n^E$-parallel nondegenerate two-form $\o_E\in \G(\La^2E)$
is fixed, then we can define a $\n$-parallel complex Riemannian metric
$\,g=\o_E\ot\o_H\,$  on $M$. 
We do not assume, in general, that the linear connection $\n$ is torsion-free.

The torsion of a linear connection belongs to
$TM\ot \La^2 T^*M$. 
Since $T^*M=E\ot H$, we have the decomposition,
\bea
TM\ot \La^2 T^*M 
&=& TM \ot \left(\La^2 E\ot S^2H \ \op\  S^2E\ot \La^2 H \right)
\nonumber\\
          &=& E^*H^* \left( \La^2 E S^2H \ \op\  S^2E \o_H \right)
\nonumber\\  &\cong& E^*\La^2 E \left( S^3 H \ \op\   \o_H H \right)
                               \ \op\  E^* S^2 E \o_H H \ ,
\la{tordecomp}\eea
where we omit the $\ot$'s and we identify $H^*$ with $H$ using
$\o_H$. 
\bd
A half-flat connection is called {\bss{admissible}} if its
torsion tensor has no component in $E^*\ot \La^2 E\ot S^3H$.  
A half-flat Grassmann manifold $(M,\n)$ is called 
{\bss admissible}  if $\n$ is {\bss admissible}. 
\ed 

We remark that if the torsion of a half-flat connection is $E$-symmetric, 
i.e.\  if it belongs to  $TM \ot S^2E  \ot \La^2 H  = TM\ot S^2E\ot\o_H$, 
then the connection is admissible.  It follows from the above decomposition 
that the torsion tensor of any
admissible connection can be written as
$$
T(e\ot h, e'\ot h')
=  T_1(e,e')\ot \o_H(h,h') h_1\ +\  T_2(e,e')\ot\o_H(h,h_2) h'\ 
      +\ T_2(e,e')\ot \o_H(h',h_2) h
$$
where $e,e'$ are sections of $E^*\,$, $\, h_1,h_2$ are fixed sections of 
$H \cong H^*$,
$T_1\in \G(E^*\ot S^2E)$ and $T_2\in \G(E^*\ot \La^2E)$. 
This shows that admissibility of the connection means that
the torsion can be represented as the sum of two tensors
linear in $\o_H$. 

\subsection{Harmonic space}
Let $M$ be a half-flat Grassmann manifold. 
We denote by $S_H$ the $\Sp{1,\bC}$-principal holomorphic
bundle over $M$
consisting of symplectic bases of $H_m^* \cong H_m \cong \bC^2\,,\ m\in M$,
$$
S_H=\{s=(h_+,h_-)\;|\;\o_H(h_+,h_-)=1 \}\, . 
$$
The bundle $S_H\rightarrow M$ is called {\bss harmonic space} \cite{GIOS}. 
A parallel (local) section
$$ m \mapsto s_m=(h_1(m),h_2(m))\in S_H$$
defines a trivialisation
$$
M \times\Sp{1,\bC} \cong S_H\ ,
$$
given by
$$
(m, \cu ) \mapsto s_m \cu
= \left(h_+ = \sum_{\a=1}^2 h_\a u_+^\a\;,\; 
        h_- = \sum_{\a=1}^2 h_\a u_-^\a
  \right)\,
,\quad \cu=  \begin{pmatrix} u_+^1 & u_-^1 \cr  u_+^2 & u_-^2\end{pmatrix}\ ;\
                  \det\ \cu = 1 \, . 
$$
We denote by  $\dpp,\dmm,\p_0$ the fundamental vector fields on $S_H$
generated by the standard generators of $\Sp{1,\bC}$,
$$
\dpp \sim  \begin{pmatrix}0&1\cr 0&0\end{pmatrix}\ ,\quad
\dmm \sim  \begin{pmatrix}0&0\cr 1&0\end{pmatrix}\ ,\quad
\p_0 \sim  \begin{pmatrix}1&0\cr 0&-1\end{pmatrix}\ . 
$$
They satisfy the relations
$$
\left[ \dpp\,,\,\dmm \right] = \p_0\ ,\quad
\left[ \p_0\,,\,\dpp \right] = 2 \dpp\ ,\quad
\left[ \p_0\,,\,\dmm \right] = -2\dmm\ . 
$$
Consider ${\rm Mat}(2,\bC )$, the vector space of two by two
matrices. The matrix coefficients $u_{\pm}^{\alpha}$ are coordinates
on this vector space. One can easily check that the vector fields
$u_+^\a  \der{u_-^\a }\,$, $u_-^\a  \der{u_+^\a }\,$ and
$\,u_+^\a \der{u_+^\a } - u_-^\a \der{u_-^\a }\,$ annihilate the function
$\ {\det\ } \cu = \e_{\b\g} u^\b_+ u^\g_-\,$, where $\e_{\b\g}$
are the matrix coefficients of the standard symplectic form of $\bC^2$. 
Therefore these vector fields are tangent to the submanifold
$\Sp{1,\bC} = \{ {\det\ } \cu = 1\} \subset{\rm Mat}(2,\bC )$. 
One can easily prove:
\bl \la{Lemma1}
In terms of the identification, $S_H \cong M\times \Sp{1,\bC}$,
the fundamental vector fields on $S_H$ generated by the standard 
generators of $\Sp{1,\bC}$ may be written:
\[
\dpp = u_+^\a  \der{u_-^\a }\ ,\quad
\dmm = u_-^\a  \der{u_+^\a }\ ,\quad
\p_0 = u_+^\a \der{u_+^\a } - u_-^\a \der{u_-^\a }\ . 
\]
\el
We say that a function $f$ on $S_H$ has {\bss charge} $c$ if
$\p_0 f = cf$. The charge measures the  difference in the degrees
of homogeneity in $u_+$ and $u_-$.

Note that any frame $(h_+,h_-)\in S_H$ defines an isomorphism
$\bC^2 \stackrel{\sim}{\ra} H^*_m$
given by $(z^1,z^2) \mapsto z^1 h_+ + z^2 h_-$. This induces an isomorphism
$$\gsp(1,\bC) = \gsp(\bC^2) \cong S^2\bC^2  \stackrel{\sim}{\ra} S^2H^*_m
        ={\rm span_\bC}  \{  h_+^2\,,\, h_-^2\,,\, h_+\vee h_-\}\ ,
$$
where we have identified $ \gsp(\bC^2)$ with  $ S^2\bC^2$ using the 
symplectic form of $\bC^2$. The generators of $\gsp(1,\bC)$ corresponding to
$\,h_+^2\,,\, -h_-^2\,,\, -h_+\vee h_-\,$ under this identification are
precisely $\,\dpp\,,\,\dmm\,,\,\p_0\,$ respectively.

\subsection{Canonical distributions on harmonic space}

Let $S_H = \{ (h_+\,,\,h_-)| h_\pm = u^\a_\pm h_\a\ ,\
(u_\pm^\a)\in \Sp{1,\bC}\}$ be the harmonic space associated
to a half-flat  Grassmann manifold $M$. Here we have fixed
a parallel symplectic frame  $(h_1,h_2)$ of $H^*$ which
defines the trivialisation $S_H = M\times \Sp{1,\bC}$ of the holomorphic
bundle $S_H$. In particular, the matrix coefficients $u^\a_\pm$ of
$\Sp{1,\bC}$ will be considered as holomorphic functions on $S_H$. 
Together with local coordinates $(x^i)$ of $M$, we obtain
a system $(x^i,u_{\pm}^{\a})$ of  local (nonhomogeneous--homogeneous)
coordinates on $S_H$. 

For any section $e\in \G(E^*)$ we define vector fields $X_\pm^e \in \gX(S_H)$
by the formula
$$ X_\pm^e \vert_{(h_+,h_-)} = \wt{e\ot h_\pm}\ ,$$
where $\wt{Y}$ stands for the horizontal lift of a tangent vector $Y$ on $M$
with respect to the connection $\n^H$. Since the frame $h_{\a}$ is parallel, 
this horizontal lift coincides with the horizontal lift with respect to 
the splitting $S_H = M \times \Sp{1,\bC}$. This shows that the vector fields
$X^e_\pm$ are tangent to $M \times \{(h_+,h_-)\}$ and hence
annihilate $u^\a_\pm\,$. 
If $h_{\pm} = u_{\pm}^{\a}h_{\a}$ then 
$X_\pm^e = u_{\pm}^{\a} \wt{X_\a^e}\,$, where $X_\a^e := e\ot h_{\a}\,$.

There exists a canonical decomposition   
$$ TS_H = T^vS_H \op \cd_+ \op \cd_- $$ 
of the (holomorphic) tangent bundle into the vertical sub-bundle $T^vS_H$ and
two (holomorphic) distributions $\cd_+$ and $\cd_-$ spanned respectively by
vector fields $X^e_+$ and $X^e_-$ associated with sections $e$ of
the  bundle $E^*$. The vertical distribution $T^vS_H$ is spanned by the
vector fields $\p_0\,, \dpp\,, \dmm\,$, which correspond to the standard 
generators of the Lie algebra $\gsp(1,\bC)$. 

\bl \la{commutatorLemma} The vector fields $X_\pm^e \in \gX(S_H)$ satisfy 
the following commutation relations:
\bea
&&\left[\p_0, X^e_\pm \right] = \pm X^e_\pm\ ,\quad
\left[\dpm, X^e_\pm \right] = 0\ ,\quad
\left[\dpm, X^e_\mp \right] = X^e_\pm 
\nonumber\\[8pt]
&&\bigl[X_+^e,X_-^{e'}\bigr]
= X_-^{\n_{\pi_*X^e_+} e'}  - X_+^{\n_{\pi_*X^{e'}_-} e }
      - \wt{T}( \pi_* X^e_+, \pi_* X^{e'}_-)
\nonumber\\[8pt]
&&\bigl[X_\pm^e,X_\pm^{e'}\bigr]
= X_\pm^{\n_{\pi_*X^e_\pm} e'} - X_\pm^{\n_{\pi_*X^{e'}_\pm} e }
         - \wt{T}( \pi_* X^e_\pm, \pi_* X^{e'}_\pm)   \ ,
\la{XX}\eea
where $T$ is the torsion of the Grassmann connection, 
$\wt{T}(X,Y) := \wt{T(X,Y)}$ denotes the horizontal lift of the vector
$T(X,Y)$ and we have used the abbreviation $\n_Xe := \n_X^Ee\,$. 
\el

\pf
The first three equations follow from 
$\,X_\pm^e = u_{\pm}^{\a}\wt{e\ot h_{\a}}\,$ and the expression for the
fundamental vector fields given in Lemma \ref{Lemma1}. To prove the last
equation we first compute the Lie bracket of two vector fields
$X= e\ot h$ and $X' = e'\ot h$ on $M$, where $h$ is parallel:
\be
\left[ X,X' \right]\ =\  \n_X X' - \n_{X'} X - T(X,X') 
\ =\ \left( \n_X e' - \n_{X'}e \right) \ot h - T(X,X')\, . 
\la{XX'}\ee
Using this we calculate the commutator
$$\arr
\left[ X^e_\pm\ ,\ X^{e'}_\pm \right] &=&
u^\a_\pm u^\b_\pm \left( \wt{\n_{X^e_\a} X^{e'}_\b} -
\wt{\n_{X^{e'}_\b} X^e_\a}
                          - \wt{T}( X^e_\a, X^{e'}_\b) \right) \\[8pt]
&=&  ( \n_{\pi_*X^e_\pm} e'\ot h_\pm )\ \wt{} - 
(\n_{\pi_*X^{e'}_\pm} e \ot h_\pm)\  \wt{}
         - u^\a_\pm u^\b_\pm \wt{T}( X^e_\a, X^{e'}_\b) \\[8pt]
&=& X_\pm^{\n_{\pi_*X^e_\pm} e' - \n_{\pi_*X^{e'}_\pm} e }
         - \wt{T}( \pi_* X^e_\pm, \pi_* X^{e'}_\pm) \ . 
\ea$$
The expression for $[X_+^e,X_-^{e'}]$ follows similarly. 
\qed

\noindent
We shall use the abbreviation
$T(X^e_\pm , X^{e'}_\pm) := \wt{T}( \pi_* X^e_\pm, \pi_* X^{e'}_\pm)$.

\bp \la{Prop2} The following conditions are equivalent:

\noindent
{~i)} For any parallel section $h \in \G(H^*)$
the distribution $E^*\ot h$ on $M$ is integrable. 

\noindent
{ii)} The distribution $\cd_{\!+}$
(associated to any parallel frame $(h_1, h_2)$)
on $S_H$ is integrable. 

\noindent
{iii)} The distribution $\cd_{\!-}$ on $S_H$ is integrable. 

\noindent
{iv)} The holomorphic
Grassmann structure is admissible, i.e.\  it has
admissible connection. 
\ep
\pf
The formula \re{XX'}, where $h$ is parallel shows that the distribution 
$E^*\ot h$ is integrable if and only if
\be
T(E^*\ot h, E^*\ot h)\subset E^*\ot h\ . 
\la{Tint}\ee
Using the decomposition \re{tordecomp}, one can check that
this condition is satisfied for all parallel sections $h$
if and only if the connection is admissible. 
This proves the equivalence of i) and iv). 
Since $\pi_* \left( X^e_+ |_{(h_+,h_-)} \right) = e\ot h_+$,
the last equation in \re{XX} shows that the distribution
$\cd_{\!+}$ is integrable if and only if \re{Tint} holds for all $h$. 
Thus i) is equivalent to ii). The equivalence of i) and iii)
is proved similarly. 
\qed

\section{Construction of half-flat connections over
half-flat Grassmann manifolds}\la{construction}

\subsection{Half-flat connections over half-flat Grassmann manifolds}
In this section we describe the {\it harmonic space method} \cite{GIOS}
for constructing half-flat connections $\n$ (Def. \ref{half-flatDef} below)
in a holomorphic vector bundle $\nu : W\ra M$ over a complex manifold $M$
with admissible half-flat holomorphic Grassmann structure.  The basic
ingredient of the construction is the lift of geometric data from $M$ to
$S_H$ via $\pi : S_H \ra M$. Let $\n$ be a holomorphic connection in a
holomorphic vector bundle $\nu : W\ra M$. Its curvature 
\be F(e\ot h_\a\,,\,e'\ot h_\b)      
 = \o_H(h_\a,h_\b) F^{(ee')} + F^{[ee']}_{\a \b }\, ,
\la{FMEqu}\ee 
where $(h_1,h_2)$ is the fixed parallel  local frame of $H^*$ and $e, e'$
are local sections of $E^*$. The curvature component  $F^{(ee')}$ is 
symmetric in $e,e'$ and $F^{[ee']}_{\a \b }$ is skew in $e,e'$
and symmetric in $\a ,\b$. 
Lifting \re{FMEqu} to $S_H$ we obtain the curvature of the pull-back
connection $\pi^*\n$ in $\pi^*\nu : \pi^*W \ra S_H$ with components,
$F(v\,, \,\cdot\, ) = 0\ ,\  \forall v\in T^vS_H $, together with
\bean
F(X^e_\pm , X^{e'}_\pm) &=& F_{\pm\pm}^{[ee']} 
:= u_\pm^{\a}u_\pm^{\b}F_{\a\b}^{[ee']} \\
F(X^e_+ , X^{e'}_-) &=&  F^{(ee')} + F_{+-}^{[ee']} := F^{(ee')} +
u_+^{\a}u_-^{\b}F_{\a \b}^{[ee']}\ .   
\eean

\bd \la{half-flatDef}
A holomorphic connection $\n$ in a holomorphic vector bundle $\nu : W\ra M$
over a complex manifold $M$ with holomorphic Grassmann structure, 
$T^*M=E\ot H$, is called {\bss half-flat} if its curvature $F$ satisfies the
equation
\be
F(e\ot h_\a\,,\, e'\ot h_\b) = \o_H(h_\a,h_\b) F^{(ee')}\, ,
\la{hfym}\ee
where $(h_1,h_2)$ is a parallel local frame of $H^*$ and
$F^{(ee')}$ is  symmetric in the local sections $e,e'$ of $E^*$. 
\ed
Note that \re{hfym} is equivalent to \re{sdgs}. 
{}From this definition it follows that for any $ h\in H^*$ we have
$
F(e\ot h\,,\, e'\ot h ) = 0\ . 
$

\bd \label{half-flatSDef} 
A connection in a holomorphic vector bundle
$W\ra S_H$ over harmonic space $S_H$ is called {\bss half-flat} if
its curvature $F$ satisfies the equations
\bea
F(X^e_+ , X^{e'}_+) &=& 0 \nonumber\\
F(X^e_+ , X^{e'}_-) &=&  F^{(ee')} \nonumber\\
F(X^e_- , X^{e'}_-) &=& 0  \nonumber\\
F(v \,,\,\cdot\ ) &=& 0 \quad,\quad \forall v\in T^vS_H\  , 
\la{hf}\eea
where $F^{(ee')}$ is  symmetric in the local sections $e,e'$ of $E^*$.
\ed

\bd \la{centralframe}
Let $\nu : W \ra M$ be a holomorphic vector bundle and $\n$ a connection
in $\pi^*\nu : \pi^*W \ra S_H$, where $\pi : S_H \ra M$. A local frame
of $\pi^*\nu$ defined on $\pi^{-1}(U)$, where $U$ is an open subset of $M$,
is called a  {\bss central frame} with respect to $\n$ if it is parallel along
the fibres of the bundle $\pi: S_H\ra M$. 
\ed
{\bf Remark:\ } If $\chi = (\chi_1,\ldots\chi_r)$ is a local frame of $\nu$
then $\pi^* \chi$ will be a central frame with respect to the pull-back
$\pi^*\n$ of any connection $\n$ in $\nu$.  The connection 1-form $A$ of 
$\pi^*\n$  with respect to the frame $\pi^* \chi$ 
satisfies $ A(v) = 0\ ,\  A(X^e_\pm) = u^\a_\pm A^e_\a\,$,
where $v$ is any vertical vector and $A^e_\a = A(\wt{X^e_\a})$ is a 
matrix-valued function on $M$. Conversely, any connection satisfying
these conditions is the pull-back of the connection over $M$ with
potential $A(X^e_\a) = A^e_\a$. 

\bp \label{pull-backProp} Let $\pi : S \ra M$ be any fibre bundle with
simply connected fibres over a simply connected manifold $M$. 
There is a natural one-to-one correspondence between gauge equivalence
classes of connections $\n^M$ in the trivial bundle  $\bC^r {\times} M$ and
gauge equivalence classes of connections $\n^S$ in $\bC^r {\times} S$
satisfying the curvature constraint $F(v,\,\cdot\,  ) = 0$ for all vertical vectors
$v$. \ep

\pf It is clear that the pull-back $\n^S = \pi^*\n^M$ to $S$
of a connection $\n^M$ defined over $M$ satisfies the curvature constraint. 
To prove the converse,
we will apply the following elementary lemma to the connection 1-form $A$
of a connection $\n$ over $N=S$. 

\bl
Let $\pi: N\ra M$ be a submersion with connected fibres
and $\a$ a $p$-form on $N$. Then
$\a$ is the pull-back $\pi^*\b$ of a $p$-form $\b$ on $M$ if and only if
the inner products $\i_v \a = \i_v d\a =0$ for all vertical tangent vectors
$v$. 
\el
Since the connection $\n^S$ is flat along the (simply connected) fibres
of $\pi$  there exists a central frame $\psi = (\psi_1,\ldots ,\psi_r)$
for $\n^S$. Let $A$ be the connection 1-form of $\n^S$ with respect to
this central frame. We then have $A(v)=0$ for any vertical vector $v$
and the  curvature condition $F(v,\,\cdot\,  ) = 0$ implies  
$ dA(v,\,\cdot\, ) = 0$. 
Now the above lemma shows
that $A$ is the pull-back of a 1-form $B$ on $M$, which defines a
connection $\n^M$ in the trivial bundle $\bC^r {\times} M$. Since any
two central frames differ by a gauge transformation which is a
matrix-valued function on $M$ the connection $\n^M$ is well defined
up to a gauge transformation. The pullback
$\pi^*\n^M$ is gauge equivalent to $\n^S$ since it has the same expression
with respect to the standard  frame of $\bC^r {\times} S$ (which is the
pull-back of the standard  frame of $\bC^r {\times} M$) as $\n^S$
with respect to the central frame $\psi$. It is clear that the
pull-backs of gauge equivalent connections over $M$ are
gauge equivalent connections over $S$. Applying a gauge transformation
to a connection $\n^S$ which has connection 1-form $A$ with respect
to a central frame $\psi$ we get a new connection $(\n^S)'$, which has the
same connection form $A$ with respect to the transformed frame $\psi'$. 
The frame $\psi'$ is therefore central with respect to $(\n^S)'$
and the two connections $\n^S$ and $(\n^S)'$ define the same
gauge equivalence class of connections over $M$. 
\qed

\bp\la{G-equiv}
Let $\nu: W = \bC^r \times M \ra M$ be a trivial vector bundle over
a complex manifold $M$ with admissible half-flat
holomorphic Grassmann structure and
$\pi^*\nu: \pi^*W = \bC^r \times S_H \ra S_H$ its pull-back to $S_H$. 
Then any half-flat connection over $S_H$ is gauge equivalent
to the pull-back
of a half-flat connection over $M$. 
\ep
\pf
It is clear that the pull-back of a half-flat connection is half-flat. 
To prove the converse, we apply Proposition \ref{pull-backProp}, by which
a half-flat connection  $\n^S$  over $S_H$ is gauge equivalent to a pull-back
connection $\pi^* \n^M$, which is necessarily half-flat. This implies that
$\n^M$ is half-flat. In fact, if the connection $\n^M$ were not half-flat then
it would have a nontrivial curvature component $F_{\a \b}^{[ee']}$ which would
imply that its pull-back $\pi^*\n^M$ has for instance a nonzero curvature
component $F_{++}^{[ee']}$. But this is impossible since $\pi^*\n^M$ is
half-flat. \qed

\bc\la{pb}
The connection 1-form $A$ of a half-flat connection over $S_H$ with respect
to a central frame $\psi$ has the form
$$ A(v)=0 \quad,\quad A(X^e_\pm) = u^\a_\pm A^e_\a\  , $$
where $v$ is any vertical vector and $A^e_\a=A(\wt{X^e_\a})$ is a
matrix-valued function on $M$. 
\ec
{\bf Remark:} This shows that the half-flat connection is completely determined
by the potential in the $\cd_{\!+}$--direction, $A(X^e_+) = u^\a_+ A^e_\a$, 
with respect to a central frame. 

\pf
This follows from Proposition \ref{G-equiv} and the remark following
Definition \ref{centralframe}. 
\qed

\subsection{The construction}
\la{constr}

In this section we construct half-flat connections in a bundle
$\nu: W\ra M$ over
a manifold $M$ with a half-flat admissible Grassmann
structure. First we define the weaker notion
of an almost half-flat connection over $S_H$ and show how to
construct all such connections from appropriate prepotentials. 
Then we show that any almost half-flat connection over $S_H$ may be
used to construct a half-flat connection on $M$. 
Since our construction is local in $M$, we shall assume that the bundles
$\pi$,
$\nu$ and $\pi^*\nu$ are trivial, i.e.\   $\pi: M \times \Sp{1,\bC}
\rightarrow M$,  $\ \nu: M \times \bC^r \ra M\,$ and
$\ \pi^*\nu: S_H \times \bC^r \ra S_H$.

\subsubsection{Construction of almost half-flat connections}

\noindent
The restriction of a half-flat connection to a leaf of the integrable
distribution ${\langle}\cd_{\!+}, \p_0{\rangle }$ is clearly flat. 
\bd
A frame $\f_1,. . . ,\f_r$ in the holomorphic vector
bundle $\pi^*\nu: \bC^r \times S_H$ which is parallel along leaves of
the integrable distribution  $\  \langle \cd_{\!+}, \p_0\rangle \ $  is called an
{\bss analytic frame}. 
\ed
With respect to an analytic frame a connection in the vector bundle
$\pi^*\nu$ has components
\bean
&\nz &=\  \p_0 \la{hfanz} \\
&\n^S_{X^e_+} &=\  X^e_+ \la{hfanp}\\ 
&\npp &=\   \dpp + \App\ := \dpp + A(\dpp )\la{hfanpp}\\
&\nmm &=\   \dmm + \Amm\ := \dmm + A(\dmm ) \la{hfanmm}\\
&\n^S_{X^e_-} &=\  X^e_- + A(X^e_-)\ . \la{hfanm}
\eean
\bd A connection $\n^S$ over $S_H$ is called {\bss almost half-flat} if its
curvature satisfies the following equations:
\bea
&& F(X^e_{+}\,,\, X^{e'}_{+}) = F(X^{e}_{+}\,,\, v ) =  0
\quad,\quad\forall v\in T^vS_H  \nonumber\\ 
&& F(\dpp\,,\,\cdot\,  ) = F(\p_0\,,\,\cdot\, )=  0\, . 
\label{almosthalf-flatEqu}\eea
\ed
In fact, these equations are not independent; for instance the Bianchi 
identity with arguments $(X_+,\dpp ,\dmm)$ together with $F(\dpp,\dmm )=
F(\dpm , X_+^e)=0$ implies the equation $F(\dpp , X^e_- ) = 0$. 

\bp \la{Propahf}
Any almost half-flat connection satisfies the following equation:
\[ F(X^e_+ ,X^{e'}_-) = F(X^{e'}_+ ,X^{e}_-)\, . \]
\ep
\pf
Using the integrability of $\cd_{\!+}$ and
$F(X_+^e, X_+^{e'} ) = F(\dmm , X_+^e) = 0$, we obtain
$$\arr
0&=& \bigl[\n^S_{\dmm}\,,\, F(X_+^e, X_+^{e'} )\bigr]\\[8pt]
&=&
\Bigl[\n^S_{\dmm}\,,\,  \Bigl[ \n^S_{X^e_+}\,,\, \n^S_{X^{e'}_+}
\Bigr] \Bigr]
- \Bigl[ \n^S_{\dmm}\,,\, \n^S_{ \left[X^e_+\,,\, X^{e'}_+ \right]}  \Bigr]
\\[10pt]
&=&
     \Bigl[ \n^S_{X^e_-}\,,\, \n^S_{X^{e'}_+} \Bigr] +
      \Bigl[ \n^S_{X^e_+}\,,\, \n^S_{X^{e'}_-} \Bigr]
                      - \n^S_{ \left[X^e_-\,,\, X^{e'}_+ \right]}
                      - \n^S_{ \left[X^e_+\,,\, X^{e'}_- \right]} \\[10pt]
&=& F(X^e_- , X^{e'}_+) - F(X^{e'}_- , X^e_+)\ . 
\ea$$
\vskip -24pt\qed

\vskip 10pt
\noindent
It follows that an almost half-flat connection is a generalisation of
a half-flat connection, satisfying only those equations
in \re{hf}, which involve curvatures with $\p_0,\dpp$ or $X^e_+$ in one
of the arguments.

\bp An almost half-flat connection is  half-flat if and only if it
satisfies $F(\dmm, X_-^e) =0$. 
\ep

\pf
By Proposition \ref{Propahf} an almost half-flat connection
is required to satisfy all the half-flatness equations \re{hf} with the
exception of
\be
F(\dmm, X_-^e) =0\quad \mbox{and} \quad F( X_-^e, X_-^{e'})=0 . 
\la{bad}\ee
The second equation here follows from the first in virtue of the
Bianchi identity with arguments $(X_+^e,X_-^{e'},\dmm)$. \qed

The following proposition shows that an almost half-flat
connection is completely determined by the potentials $A_{++}$ and
$A_{--}$ with respect to an analytic frame. 
\bp \label{ahfProp1}
Let $\n^S$ be an almost half-flat connection in the vector bundle
$\pi^*\nu: \bC^r {\times} S_H \ra S_H$ with potentials
$A_{++}$, $A_{--}$ and $A(X_-^e)$ in an analytic frame. Then:

\noindent
{\it (i)\ } The potential $A_{++}$ is analytic and has charge $+2$, i.e. 
\be \label{prepEqu}
X^e_+ \App = 0\quad ,\quad  \p_0 \App = 2 \App\ . 
\ee
{\it (ii)\ }  The potential $A_{--}$ satisfies
\be
\dpp \Amm - \dmm \App + [ \App , \Amm ] = 0 \quad,\quad
\p_0 \Amm =-2 \Amm\ . 
\la{amm}\ee
{\it (iii)\ } The potential $A(X_-^e)$  
is determined by $A_{--}$ and has charge $-1$:
\be \label{amEqu}
   A(X^e_-) =  - X^e_+ \Amm\quad ,\quad  \p_0 A(X^e_-) =- A(X^e_-)\ . 
\ee
Conversely, any matrix-valued potentials  $A_{++}$, $A_{--}$ and $A(X_-^e)$
satisfying \re{prepEqu}, \re{amm} and  \re{amEqu} define an
almost half-flat connection. 
\ep

\pf
(i) The curvature constraints
$F(X^e_+\,,\,\dpp)=0\ ,\ F(\p_0\,,\,\dpp)=0$, in an analytic frame,
take the form \re{prepEqu}. \\
(ii) The further almost half-flatness conditions,
$F(\dpp\,,\,\dmm) = F(\p_0\,,\,\dmm)=0$ give
equations \re{amm} for the potential $\Amm\,$. \\
(iii) Having obtained $\,\Amm\,$, we can find $\,A(X^e_-)\,$ from
the equations $\,F(X^e_+\,,\,\dmm) = F(\p_0\,,\,X^e_-) = 0\,$, which
take the form
\be
   X^e_+ \Amm = A([X_+^e,\dmm]) = - A(X^e_-) \quad,\quad
\p_0 A(X^e_-) =- A(X^e_-)\ .
\la{am}\ee
The second equation follows from the first. 
\qed

We can now write an algorithm for the construction of all
almost half-flat connections: 
\bt\la{algo}
Let $\App$ be an analytic prepotential, i.e.\  a matrix-valued
function on a domain $U = \pi^{-1}(V) \subset S_H$, where $V\subset M$ is
a simply connected domain, satisfying \re{prepEqu}. Let $\F$ be an
invertible matrix-valued function on $U$ which satisfies the
equations
\be \label{b} \dpp \Phi = -A_{++}\Phi \, ,\quad \p_0 \Phi = 0\, .
\ee
It always exists.
The pair $(\App , \F)$ determines an almost  half-flat connection
$\n^S = \n^{(\App , \F)}$. Its potentials with respect to an analytic frame
are given by $\App$, $\Amm = -(\dmm \F) \F^{-1}$ and
$A(X^e_-) = - X^e_+ \Amm\, $. Conversely, any almost  half-flat connection
is of this form. 
\et

\pf We consider the connection defined by $\App$ and $A(\p_0) = 0$
along an orbit $sB$ of the Borel subgroup of SL($2,\bC$),
\[ B = \left\{ \left. 
                \begin{pmatrix}t_0 & t_1\cr
                        0 & t_0^{-1}\end{pmatrix}\right| \, t_0\in \bC^*, \,
t_1 \in \bC \right\}
\cong \bC^*\times \bC \quad \mbox{(diffeomorphic)}\ .
\]
It is flat since the second equation of \re{prepEqu}
is equivalent to $F(\p_0\,,\,\dpp){=}0$ (vanishing of the curvature along
$sB$). Moreover it has trivial holonomy since
the fundamental group of $B \cong \bC^*{\times} \bC$ coincides with the
fundamental group of the $\bC^*$-factor and the potential is zero in
direction of $\p_0$ which is tangent to $\bC^*$. 
An invertible solution to the system \re{b} exits and
defines a parallel frame $\F$ with respect to the
flat connection with trivial holonomy defined along each orbit of the
Borel group. Since the space of Borel orbits in $U$ is diffeomorphic
to $V \times \bC P^1$ and is therefore simply connected, a solution $\F$
exists on the domain $U$. 
Now, given any such solution of \re{b}, we define
$A_{--} := -(\dmm \F ) \F^{-1}$. This solves \re{amm}, since
$F(\dpm,\p_0 ) = F(\dpp ,\dmm ) = 0$ is the integrability condition
for the system $\dpm \F = -\Apm \F$, $\p_0 \F = 0$. Finally, we define
$A(X^e_-) := - X^e_+ \Amm\,$, obtaining an almost half-flat connection
by Proposition \ref{ahfProp1}. Now the converse statement follows also
from Proposition \ref{ahfProp1}. 
\qed

\subsubsection{Transformation to the central frame}

Since an almost half-flat connection $\n = \n^S$ is flat in vertical
directions, it admits a central frame $\psi$. 
The following lemma shows that the solution $\F$ of equation \re{b} gives a
gauge transformation from an analytic frame $\f$ to a central frame 
$\psi = \f\F$  for the almost half-flat connection $\n^{(A_{++},\F)}$. 
\bl\la{centralLem}
Let $\n = \n^{(A_{++},\F)}$ be the almost  half-flat
connection associated to the analytic prepotential $A_{++}$ with respect to
the analytic frame $\f$ and an invertible solution $\F$ of $\re{b}$. 
Then the frame $\psi := \f\F$ is a central frame for the connection $\n$,
i.e.\  the potentials $C(\dpm )$ and $C(\p_0)$ with respect to that frame
vanish. 
\el
\pf The result follows from the pure gauge form of $A(\dpm )$ and $A(\p_0)$
and the transformation law \re{gt} for potentials. 
\qed

\noindent
With respect to the central frame $\j$, the almost half flat
connection constructed above takes the form:
$$\left\{ \arr
&\n^S_{X^e_+} &=\ X^e_+ + C({X^e_+})\
                  =\ X^e_+ + \F^{-1}\ X^e_+\ \F  \\[6pt]
&\n^S_{X^e_-}&=\  X^e_- + C({X^e_-})\ =\  X^e_- +  \F^{-1} X^e_- \F
                     + \F^{-1} X^e_+\left(\dmm\F\F^{-1}\right)\F  \\[6pt]
&\npp &=\   \dpp  \quad,\quad
\nmm\ =\   \dmm  \quad,\quad
\nz\ =\  \p_0\ .
\ea \right. $$
Moreover, the equations $F(\dpp\,,\,X^e_+) = F(\p_0 \,,\,X^e_+) =0$ 
imply that the potential $C({X^e_+})$ satisfies the equations 
\be
\dpp C({X^e_+}) = 0\ ,\quad \p_0 C({X^e_+}) = C({X^e_+})\ .
\la{hfCp}\ee

\subsubsection{Construction of half-flat connections}

We assume now that the analytic prepotential $\App$ is defined globally
along the fibres of $\pi: S_H \ra M$. Then, restricting $M$ to an appropriate
domain, we may assume that  $\App$ is defined globally on $S_H$. The previous
construction then provides an almost half-flat connection over $S_H$. 
Using this connection, we may construct a half-flat connection on $M$. 
The crucial point is the following:
\bp\la{linC}
The potential $C({X^e_+})$ of an almost half-flat connection $\n$
with respect to a central frame  is linear in $u_+^\a\,$, namely,
\be
C({X^e_+}) = u_+^\a C(\wt{X^e_\a}) =: u_+^\a C^e_\a\ ,  
\ee
where $(x^i, u_\pm^\a)$ are the local coordinates 
associated with the trivialisation
$S_H = M\times \Sp{1,\bC}$ and $C^e_\a = C^e_\a(x^i)$ is a matrix-valued 
function on $M$. 
\ep
\pf
Due to equations \re{hfCp}, the result follows from:
\bl\la{lemf+}
i) If a holomorphic function $f_+\,$, defined on some 
domain 
\[ U \subset \{ u_+^2 \neq 0\} \subset \Sp{1,\bC}= \left\{\cu=
 \begin{pmatrix} u_+^1 & u_-^1 \cr  u_+^2 & u_-^2\end{pmatrix}\ ;\
                  \det \cu = 1 \right\} \]
satisfies
\be
\dpp f_+ =0\ ,\quad  \p_0 f_+ = f_+\ ,
\la{syst}\ee
then $\ f_+ = u_+^\a  f_\a(\fr{u^1_+}{u^2_+})$. 
Here $f_\a(\fr{u^1_+}{u^2_+})$ are holomorphic  functions on
$U$ invariant under the right action of the Lie algebra of 
upper-triangular matrices. 

\noindent
ii) Moreover, if the function $f_+$ is globally defined, then it
is linear in $u^\a_+$, i.e.\    $\ f_+ = u_+^\a  f_\a, \quad f_\a =$const. 
\el
\pf
i) One can immediately check that $\,f_+ = u_+^\a  f_\a(\fr{u^1_+}{u^2_+})\,$
is a solution of \re{syst}. 
We note that the quotient of any two solutions of \re{syst} is
a solution of the corresponding homogeneous system,
\be
\dpp f =0\ ,\quad  \p_0 f = 0\ . 
\la{syst_hom}\ee
It is sufficient to check that any solution of \re{syst_hom} is
a function of $u^1_+/u^2_+\,$. To prove this we use the local
factorisation of $\Sp{1,\bC}$ into the product of a Borel subgroup $\cb$
and a nilpotent subgroup as follows:
$$
 \begin{pmatrix} u_+^1 & u_-^1 \cr  u_+^2 & u_-^2\end{pmatrix} =
  \begin{pmatrix}1&0\cr c^{-1}&1\end{pmatrix}
  \begin{pmatrix}a&b \cr 0&a^{-1}\end{pmatrix}\quad ,\quad
\cb= \left\{  \begin{pmatrix}a&b \cr 0& a^{-1}\end{pmatrix} \right\}\ .
$$
Then $\,c=u^1_+/u^2_+\,$ and $\,\p_0,\dpp\,$ are generators of the right
action of $\cb$. This implies that the solutions of \re{syst_hom} are
precisely the local functions on $\Sp{1,\bC}$ invariant under the right
action of $\cb$.  In terms of the local coordinate system $(a,b,c)$ on 
$\Sp{1,\bC}$ such functions are functions of $c=u^1_+/u^2_+\,$ alone. 

\noindent
ii) The restriction $V|_{\Sp{1}}$ to $\Sp{1}$ of any
irreducible $\Sp{1,\bC}$-module $V$ of holomorphic functions
is a (finite dimensional) irreducible $\Sp{1}$-module  of smooth
functions on $\Sp{1}$. 
The condition \re{syst} shows that $f_+$ is a highest weight vector
with weight $+1$. Hence $f_+$ generates a 2-dimensional submodule 
$\langle f_+\rangle  = {\rm span}\{ f_+ , f_-:=\dmm f_+\}$ of holomorphic functions. 
It remains to show that any two dimensional
module of holomorphic functions on $\Sp{1,\bC}$ is spanned by
linear functions. We know two such modules, generated by
the highest weight vectors $u_+^1$ and $u_+^2$ respectively. 
On the other hand, by the Peter-Weyl Theorem the multiplicity
of the two-dimensional irreducible representation of $\Sp{1}$
in $L^2(\Sp{1})$ is 2. 
\qed

Using Proposition \ref{linC}, with respect to a central frame, we can write
$\n_{X^e_+} = X^e_+ + u^\a_+ C^e_\a $ where the coefficients
$C^e_\a = C^e_\a(x^i)$
are matrix valued functions of coordinates $x^i$ on $M$. Using them
we define a new connection in the trivial bundle $\bC^r \times S_H$ over
$S_H$ by
$$ \left\{ \arr
\wh\n_{X^e_\pm} &=& X^e_\pm + u^\a_\pm C^e_\a \\[4pt]
\wh\n_{\dpm} &=&  \dpm \quad,\quad  \wh\n_{\p_0}\ =\ \p_0\ .
\ea \right. $$
Our main result now follows:
\bt \label{hatThm} Let $M$ be
a complex manifold with a half-flat admissible Grassmann
structure. Let $A_{++}$ be an analytic prepotential, i.e.\  a solution
of \re{prepEqu}, and $\F$ an invertible solution of \re{b}. 
Then the connection $\wh\n = \wh{\n}^{(A_{++},\Phi )}$ constructed from the 
data $(A_{++},\F )$
is a half-flat connection in the trivial vector bundle
$\bC^r \times S_H \ra S_H$
and it is the pull-back of the following half-flat connection
$\n^M$ in the bundle $\bC^r \times M\ra M$:
\be
\n^M_{X^e_\a} = X^e_\a  + C^e_\a\ .
\ee
Conversely, any half-flat connection over $S$ (or $M$)
is gauge equivalent to one
obtained from the above construction. 
\et
\pf
The remark after Definition \ref{centralframe} shows that the connection 
$\wh\n$ is the pull-back of
the connection $\n^M$. It suffices now to show that $\n^M$ is half-flat. 
Note that the connections $\n$ and $\wh\n$ coincide in the direction of
$X^e_+$. Hence, using $ C^e_+ := u^\a_+ C^e_\a$, we have
$$\arr
0 &=&  F^{\n}(X^e_+ , X^{e'}_+)\
   =\ F^{\wh\n}(X^e_+ , X^{e'}_+) \\[6pt]
     &=&  X^e_+ C^{e'}_+ - X^{e'}_+ C^e_+ + \left[ C^e_+ \,,\, C^{e'}_+
\right]
                   - C({\left[ X^e_+ \,,\, X^{e'}_+ \right]}) \\[6pt]
      &=& u^\a_+ u^\b_+ \left( X^e_\a C^{e'}_\b - X^{e'}_\b C^e_\a +
\left[ C^e_\a \,,\, C^{e'}_\b \right]
                   - C({\left[ X^e_\a \,,\, X^{e'}_\b \right]}) \right)
\\[6pt]
      &=& u^\a_+ u^\b_+  F^{\n^M}( X^e_\a \,,\, X^{e'}_\b )\ ,
\ea$$
since $X^e_+ u^\b_+ = 0$. 
This shows that the curvature $F^{\n^M}( X^e_\a \,,\, X^{e'}_\b )$ is
skewsymmetric in $\a,\b$, i.e.\  it belongs to $\La^2H\ot S^2E\ot \End W$. 
In other words, the connection $\n^M$ is half-flat. 

Conversely, let $\n^S$ be a half-flat connection over $S_H$. 
By Proposition \ref{G-equiv} we may assume that it is a pull-back of
a half-flat connection $\n^M$ over $M$. Since the restriction of $\n^S$ to
the leaves of $\langle \cd_{\!+}, \p_0\rangle $ is flat, there exists an analytic frame
(i.e.\ a frame such that $A(X_+^e) = A(\p_0) = 0$, in which the potential
$A(\dpp )$ satisfies the equations \re{prepEqu}). Since $\n^S$ is flat along
the (simply-connected) fibres, there exists an invertible solution $\F$ to 
the system
\[ \dpp \F + A_{++}\F\ =\ \dmm \F + A_{--}\F\ =\ \p_0 \F\ =\ 0 \, . 
\]
This shows that $\n^S$ is gauge-equivalent to the almost half-flat connection
$\n^{(\App,\F)} = \wh{\n}^{(\App,\F)}$. 
\qed

\subsection{Application to hyper-K\"ahler manifolds with admissible
torsion}
The above construction can be applied to the complexification
of \hk\ manifolds. Recall that any \hk\ manifold admits a
(locally defined)  Grassmann structure $T^{*\bC} M = E\ot H$, 
such that the Levi-Civita connection on the cotangent bundle 
$\n = \n^E \ot {\rm Id} + {\rm Id} \ot \n^H$ is half-flat, 
i.e.\  the connection $\n^H$ is flat. Since the \hk\ metric
is Ricci flat, hence analytic, we may, using
analytic continuation, extend the manifold $M$ to a complex
manifold $M^\bC$ with holomorphic extension of the \hk\ structure,
in particular, we have a holomorphic Ricci flat metric on $M^{\bC}$
with holonomy in $\Sp{n,\bC}$ and half-flat Grassmann structure. 
This Grassmann structure is admissible since the Levi-Civita
connection on $M^{\bC}$ has no torsion. Hence we can apply the
harmonic space method to construct half-flat connections on holomorphic
vector bundles $W \ra M^{\bC}$. The complex version of Proposition
\ref{l1sdProp} shows that such connections are Yang-Mills connections. 
More generally, the method of construction of half-flat connections
extends to real analytic
(possibly indefinite) \hk\ manifolds
with admissible torsion, i.e.\  with torsion which has
zero component in $S^3H \ot E^*\ot \La^2 E$. 
A {\bss \hk\ manifold with admissible torsion} is defined as a
pseudo-Riemannian manifold $(M,g)$ with a linear metric connection $\n$
with holonomy in $\Sp{k,l}$ which has admissible torsion. 
As in the (torsionfree) \hk\ case  there exists a
parallel $4$-form given by $\O =\sum_\a \o_\a \wedge \o_\a\,$,
$\o_\a := g(J_\a \cdot ,\cdot)$, and half-flat connections are
characterised  as connections with curvature in
$V_{\l_1} \ot \End W$, where $V_{\l_1}$ is the $\l_1$-eigenspace
of the endomorphism $B_{\O}$ associated to  $\O$. If the form
$\O$ is co-closed then any half-flat connection will be
$(\O ,\l_1)$-self-dual and thus a Yang-Mills connection. 
We remark that co-closedness of $\O$ is equivalent to a linear
$\Sp{k,l}$-invariant condition on the torsion.

\section{Generalisation to higher-spin Grassmann manifolds}

\subsection{Higher-spin Grassmann structures}

The construction discussed in the previous section is in fact the $m{=}1$
specialisation of a more general construction of connections on
{\it spin~$\fr{m}{2}$ Grassmann manifolds}, which we discuss in this section. 
These manifolds were considered in \cite{AG}. 
\bd
A {\,\bss spin $\fr{m}{2}\,$ Grassmann structure\,} on a (complex)
manifold $M$ is a holomorphic Grassmann structure of the form
$T^*M \cong E\ot F =  E\ot S^m H$, with a holomorphic Grassmann connection
$\n = \n^E \ot {\rm Id} + {\rm Id} \ot \n^F$, where $H$ is a rank 2
holomorphic vector bundle over $M$ with holomorphic symplectic connection 
$\n^H$ and symplectic form $\o_H\in\G(\La^2 H)$, and $\n^F$ is the connection
in $F=S^mH$ induced by $\n^H$. 
$M$ is called {\bss half-flat} if the connection $\n^F$ is flat. 
\ed
The bundle $S^m H$ is associated with the
spin~$\fr{m}{2}$ representation of the group $\Sp{1,\bC}$. 
Any frame $(h_1,h_2)$ for $H^*$ defines a frame for $S^m H^*$,
$(h_A := h_{\a_1}h_{\a_2}\dotsm h_{\a_m})$, where the multi-index
$A:=\a_1\a_2\ldots\a_m\ ,\ \a_i=1,2$. 
The $\n^H$-parallel symplectic form $\o_H$ on $H^*$ induces a bilinear form
$\o_H^m$ on  $F^*=S^m H^*$ given by
$$
\o_H^m(h_A,h_B) := \underset{A}{\gS}\, \underset{B}{\gS}\;
\o_H(h_{\a_1},h_{\b_1}) \o_H(h_{\a_2},h_{\b_2})\dotsm
\o_H(h_{\a_m},h_{\b_m})\ , $$
where $\gS_A$ denotes the sum over all permutations of the $\a$'s. 
This form is skew-symmetric if $m$ is odd and symmetric if $m$ is even. 
To any section $e\in \G(E^*)$ and multi-index $A$ we associate
the vector field $X^e_A := e\ot h_A$ on $M$. 

The construction of half-flat connections described in section \ref{constr}
may be adapted to obtain certain `partially flat' connections in vector
bundles $W\ra M$, provided that the torsion of $\n$ obeys
certain admissibility conditions. 

\bd
Let $(M,\n )$ be a half-flat spin~$\fr{m}{2}$ Grassmann manifold. 
For any section $e\in \G(E^*)$ we define vector fields
$$\arr  X^e_{(m-2i)+} &:=&
u^{\a_1}_- \dotsm u^{\a_i}_- u^{\a_{i+1}}_+ \dotsm u^{\a_m}_+ X^e_A\
\quad \mbox{if} \quad m-2i \ge 0\\[6pt]
X^e_{(2i-m)-} &:=&
u^{\a_1}_- \dotsm u^{\a_i}_- u^{\a_{i+1}}_+ \dotsm u^{\a_m}_+ X^e_A\
\quad \mbox{if} \quad m-2i <0
\ea$$
on the principal bundle $S_H$ of symplectic frames in $H$. The
distribution spanned by these vector fields is denoted by
$\,\cd_{k+} :=\, \langle X^E_{k+}\rangle\,$ 
for $k{\ge}0\,, k\equiv m\!\mod 2$ and 
$\, \cd_{k-} :=\, \langle X^E_{k-}\rangle\,$ 
for $k{>}0\,, k\equiv m \!\mod 2$. We define also 
$$
\cd^k_{(\pm)} :=  \bigoplus_{i= 0}^k \cd_{(m-2i)\pm}\ . 
$$
The Grassmann connection $\n$ is called {\bss k--admissible} 
if it preserves the distribution
$\cd^k_{(\pm)}$, i.e. 
\be
T\left(\cd^k_{(\pm)},\cd^k_{(\pm)} \right)\subset \cd^k_{(\pm)}\ . 
\la{Dkadmiss}
\ee The Grassmann manifold $(M,\n)$ is
called {\bss k--admissible} if the connection $\n$ is
$k$--admissible. 
\ed For small $m$ we shall write $X^e_0$, $X^e_+$,
$X^e_-$, $X^e_{++}$ etc. \ instead of $X^e_{0+}$, $X^e_{1+}$,
$X^e_{1-}$, $X^e_{2+}$ etc. 

\bp Let $(M,\n )$ be a half-flat
spin~$\fr{m}{2}$ Grassmann manifold. Then the distribution
$\cd^k_{(\pm)}$ is integrable if and only if the torsion of the
Grassmann connection $\n$ satisfies equation \re{Dkadmiss}. 
\ep
The proof is similar to that of Proposition \ref{Prop2}.

\subsection{Partially flat connections over higher-spin
Grassmann mani\-folds}

Let $(M,\n )$ be a half-flat spin~$\fr{m}{2}$ Grassmann manifold and
$\nu: W\ra M$ a holomorphic vector bundle. Since our constructions are
local we will assume that $W$ is trivial. In the higher spin $(m>1)$
case, there exists, as a natural generalisation of the notion of a
half-flat connection, the more refined notion of a $k$-partially
flat connection in $\nu$. The space of two-forms $\La^2 T^*M$ has
the following decomposition into $GL(E)\ot \Sp{1,\bC}$-submodules:
$$ \La^2 T^*M = \La^2\bigl( E\ot S^mH\bigr)
                 = \La^2E \ot S^2 S^mH \ \op\  S^2E  \ot  \La^2S^mH\, ,
$$
where
\bean 
S^2 S^mH &=& S^{2m}H\ \op\ \o_H^2 S^{2m-4}H\ \op\ \dotsb  \op
\o_H^{2[\frac{m}{2}]} S^{2m -4[\frac{m}{2}]}H 
\\[2pt]
\La^2S^mH &=&  \o_H S^{2m-2}H\ \op\  \o_H^3 S^{2m-6}H \op \dotsb
\op \o_H^{2[\frac{m}{2}]+1} S^{2m - 4[\frac{m}{2}]-2}H\, . 
\eean
Here we use the convention that $S^lH = 0$ if $l< 0$. 

Let $\n$ be a connection in the vector bundle $W\ra M$. 
Its curvature has the following decomposition, corresponding to
the above decomposition of $\La^2 T^*M$ into irreducible
$GL(E){\cdot}\Sp{1,\bC}$-submodules:
\bea
F(X^e_A, X^{e'}_B) 
&=& \underset{A}{\gS}\, \underset{B}{\gS}\; \sum_{k=0}^{\left[m/2\right]}
        \Bigl(
        \o_H(h_{\a_1},h_{\b_1})\dotsm \o_H(h_{\a_{2k}},h_{\b_{2k}})
   \stackrel{(2k)}{F}{}^{[ee']}_{\a_{2k+1}\ldots\a_m\b_{2k+1}\ldots\b_m} 
\bigr. \nonumber\\[2pt] &&\  + \bigl. 
        \o_H(h_{\a_1},h_{\b_1})\dotsm \o_H(h_{\a_{2k+1}},h_{\b_{2k+1}})
   \stackrel{(2k+1)}{F}{}^{(ee')}_{\a_{2k+2}\ldots\a_m\b_{2k+2}\ldots\b_m}
\Bigr),
\la{curv_decomp}\eea
where the tensors $\stackrel{(2k)}{F}\in \G(\La^2 E \ot S^{2m-4k}H)$
and  $\stackrel{(2k+1)}{F}\in \G(S^2 E \ot S^{2m-4k-2}H)$. 

We note that half-flat connections are those which satisfy the
conditions
\be
\stackrel{(2i)}{F} = 0\ ,\quad {\rm for\ all}\ i\in \bN . 
\ee
For $m{>}1$ these conditions are not suitable for application of the
harmonic space method. However, the following more refined restrictions on 
the curvature are amenable to the method (c.f.\ \cite{DN1}). 
\bd
A connection $\n$ in the vector bundle $\nu: W\ra M$ is called
{\bss k-partially flat} if $\stackrel{(i)}{F}=0$ for all $i\le 2k$. 
Here $0\le k\le [\fr{m+2}{2}]$. 
\ed
Clearly, $[\fr{m+2}{2}]$-partially flat connections are simply flat
connections. We note that for $m{=}1\,$, 0-partially flat connections
are precisely half-flat connections. 
For general odd $m{=}2p{+}1\,$, 0-partially flat connections in a
vector bundle $\nu$ over flat spaces with spin~$\fr{m}{2}$
Grassmann structure were considered by Ward \cite{W}. He chose $E$ 
to be a rank 2 flat bundle and showed that 0-partially flat
connections, for $m{>}1$, do not correspond to Yang-Mills
connections. Therefore, in our more general setting, we clearly cannot
expect 0-partially flat connections to satisfy the Yang-Mills equations
for $m{>}1$. On the other hand, the penultimate case, $k = [\fr{m}{2}]$,
is particularly interesting for odd $m$:
\bt\la{pfymTh} Let $M$ be a  half-flat
spin~$\fr{m}{2}$ Grassmann manifold $M$. If $m$ is odd and the
vector bundle $E^*\ra M$ admits a $\,\n^E$-parallel symplectic form
$\o_E$, then $M$ has canonical $\,\Sp{E}{\cdot}\Sp{H}$-invariant metric
$g = \o_E \otimes \o_H^m$ and $4$-form $\,\Omega\neq 0$. 
If $\,\O\,$ is co-closed with respect to the metric $g$, then any
$\fr{m-1}{2}$-partially flat connection $\n$ in a vector bundle
$W$ over $M$ is $(\O , \l )$-self-dual and hence it is a Yang-Mills connection. 
\et

\pf 
To describe $\O$ we use the following notation: $e_a$ is a basis of $E^*$,
$h_\a$ is a basis of $H^*$, $h_A$ is the corresponding basis of $S^mH^*$ and
$X_{aA} := e_a\ot h_A$ the corresponding basis of $TM = E^*\ot S^mH^*$. 
With respect to these bases, the skew symmetric forms $\o_E$, $\o_H$ and
$\o_H^m$ are represented by the matrices $\o_{ab}$, $\o_{\a \b }$ and
$\o_{AB}$ respectively. We define $\O$ by
\[ \O := \sum \o_{ab}\,\o_{cd}\,\o_{AC}\,\o_{BD}\,
              X^{aA}\wedge X^{bB}\wedge X^{cC}\wedge X^{dD} \, ,\]
where $X^{aA}$ is the basis dual to $X_{aA}$. This form is
obviously $\Sp{E}{\cdot} \Sp{H}$-invariant since we used only $\o_E$ and
$\o_H$ in the definition. One can easily check that $\O \neq 0$. 
The connection $\n$ is $\fr{m-1}{2}$-partially flat if and only if its
curvature $F$ belongs to the space
\[ S^2E {\otimes} \o_H^m {\otimes} \End W\
\subset\ S^2E {\otimes} \wedge^2S^mH  {\otimes} \End W\
\subset\ \wedge^2 \left(E{\ot} S^mH\right){\ot}\End W  \, . \]
Here we use the decomposition
\[ \wedge^2S^mH = \o_H S^{2m -2}H\ \oplus\ \o_H^3 S^{2m - 6}H
\ \oplus\  \ldots \ \oplus\   \bC \o_H^m \, . \]
The $\Sp{E}{\cdot}\Sp{H}$-submodule $S^2E \otimes \o_H^m \subset
\wedge^2T^*M$
is irreducible. Therefore it is contained in an eigenspace
$V_\l$ of the $\Sp{E}{\cdot} \Sp{H}$-invariant operator
$B_\O : \wedge^2T^*M \ra \wedge^2T^*M$. It remains to check that $\l \neq 0$. 
By Lemma \ref{contractionLemma} it suffices to compute the contraction
$K = K_{cCdD}X^{cC}X^{dD}$ of
a tensor $S = S_{ab}\o_{AB}e^ae^b\otimes h^Ah^B$ in $S^2E \otimes \o_H^m$
with $\O$:
\begin{eqnarray*} -K_{cCdD} &=& S^{ab}\o^{AB} (\o_{ab}\o_{cd}\o_{AC}\o_{BD} +
\o_{ac}\o_{db}\o_{AD}\o_{CB} + \o_{ad}\o_{bc}\o_{AB}\o_{DC}
-\o_{ba}\o_{cd}\o_{BC}\o_{AD}\\
& & -\o_{bc}\o_{da}\o_{BD}\o_{CA} -\o_{bd}\o_{ac}\o_{BA}\o_{DC}
 -\o_{ca}\o_{bd}\o_{CB}\o_{AD} -\o_{cb}\o_{da}\o_{CD}\o_{BA}\\
& & -\o_{cd}\o_{ab}\o_{CA}\o_{DB}+ \o_{da}\o_{bc}\o_{DB}\o_{AC} +
\o_{db}\o_{ca}\o_{DC}\o_{BA} + \o_{dc}\o_{ab}\o_{DA}\o_{CB})\\
&=& 4(m+1)S_{cd}\o_{CD} \, . 
\end{eqnarray*}
Here $\o^{ab}$ and $\o^{AB}$ denote the inverses of $\o_{ab}$ and $\o_{AB}$
and $S^{ab} = \o^{aa'}\o^{bb'}S_{a'b'}$. We have used that  $S^{ab}\o^{AB}$
is skew symmetric under interchange of $aA$ with $bB$ and that
$\o^{AB}\o_{AB}{=}-(m+1)$. The above calculation shows that 
$\l{=}-4(m+1)\neq 0$ and hence any $\fr{m-1}{2}$-partially flat connection 
is $(\O , \l)$-self-dual and is a Yang-Mills connection by Theorem \ref{YMThm}. 
\qed

\noindent
The analogous result does not hold if $m$ is even. 
\bp If  $m$ is even and the
vector bundle $E^*{\ra} M$ admits a $\n^E$-parallel metric $\g_E$, then
$M$ has canonical $\SO{E}{\cdot}\Sp{H}$-invariant metric
$g = \g_E \otimes \o_H^m$ and  $4$-form $\Omega\neq 0$. 
\ep
\pf  Analogously to the case of $m$ odd, we can define $\O$ by
\[ \O := \sum \g_{ab}\,\g_{cd}\,\o_{AC}\,\o_{BD}\,
          X^{aA}\wedge X^{bB}\wedge X^{cC}\wedge X^{dD} \, . \]
Here $\g_{ab} = \g_E(e_a,e_b)$ and we recall that for even
$m$ the bilinear form $\o_H^m$ is symmetric:
$\o_{AB} = \o_{BA}$. 
\qed

\noindent
For even $m$, a connection $\n$  in a vector bundle
$W$ over $M$ is $\fr{m}{2}$-partially  flat if and
only if its curvature $F$ belongs to the space
\bean
&& \Bigl( \wedge^2E {\ot} \o_H^m\ \op\ 
S^2E {\ot} S^2H{\ot} \o_H^{m-1}\Bigr) {\ot} \End  W \\
&& \subset \Bigl( \wedge^2E {\ot} S^2S^mH\ \oplus\ 
S^2E {\ot} \wedge^2S^mH\Bigr) {\ot} \End  W 
= \wedge^2 \left( E{\ot} S^mH \right) {\ot} \End  W . 
\eean
The $\SO{E}{\cdot} \Sp{H}$-submodule
$\wedge^2E \otimes \o_H^m \op S^2E \otimes S^2H\ot
\o_H^{m-1} \subset  \wedge^2T^*M$ is not irreducible, so
unlike the odd $m$ case we cannot conclude that it is
contained in an eigenspace
$V_\l$ of the $\SO{E}{\cdot} \Sp{H}$-invariant operator
$B_\O : \wedge^2T^*M \ra \wedge^2T^*M$. In fact, examples are known
(see Appendix B of \cite{DN2}) where $B_\O$ has different eigenvalues on each
irreducible summand of $\wedge^2T^*M$. Therefore in the case of even $m$ we
cannot expect that $\fr{m}{2}$-partial flatness implies the Yang-Mills
equations.

\subsection{Construction of partially flat connections over higher spin
Grassmann mani\-folds}
Now we generalise the construction of half-flat connections over admissible
half-flat Grassmann manifolds to the case of $k$--partially flat
connections over $k$--admissible higher spin Grassmann manifolds $M$. 
The natural extension of the harmonic construction given in section
\ref{constr} yields $k$-partially flat connections in the vector bundle
$\nu$ over the $k$--admissible spin~$\fr{m}{2}$ Grassmann manifold $M$. 
Again, we lift the
geometric data from $M$ to $S_H$ via the projection $\pi: S_H\ra M$. 
The pull back  $\pi^* \n$ of a $k$--partially flat connection $\n$
in the trivial  vector bundle $\nu: W = \bC^r \times M \ra M$ is
a connection in the vector bundle $ \pi^*\nu: \pi^* W\ra S_H$ which
satisfies equations defining the notion of a $k$--partially flat gauge
connection on $S_H$. One can also define the weaker notion of an almost
$k$--partially flat connection in $ \pi^*\nu: \pi^* W\ra S_H$. The latter
may be constructed from a prepotential and it affords the construction
of a $k$-partially flat connection in the bundle $W\ra M$. 
To simplify our exposition we explain the construction in the $m{=}3$
case. Here the decomposition \re{curv_decomp} of the curvature
tensor takes the
form
\bea 
     F(X^e_{\a_1\a_2\a_3}, X^{e'}_{\b_1\b_2\b_3})
& =& \underset{A}{\gS}\, \underset{B}{\gS}\; 
        \left( \stackrel{(0)}{F}{}^{[ee']}_{\a_1\a_2\a_3\b_1\b_2\b_3}
      + \o_H(h_{\a_1},h_{\b_1})\stackrel{(1)}{F}{}^{(ee')}_{\a_2\a_3\b_2\b_3}
\right. \nonumber\\[4pt] &&\qquad
+ \o_H(h_{\a_1},h_{\b_1})\o_H(h_{\a_2},h_{\b_2})
            \stackrel{(2)}{F}{}^{[ee']}_{\a_3\b_3}
\nonumber\\[4pt] &&\qquad \left. 
     + \o_H(h_{\a_1},h_{\b_1}) \o_H(h_{\a_2},h_{\b_2})\o_H(h_{\a_3},h_{\b_3})
       \stackrel{(3)}{F}{}^{(ee')} \right).
\label{3curv_decomp}\eea
In this case we have two nontrivial notions of partial flatness:
\bea
\text{0--partial flatness}:&&   \stackrel{(0)}{F}=0
\la{0pf}\\[10pt]
\text{1--partial flatness}:&&
 \stackrel{(0)}{F}= \stackrel{(1)}{F}= \stackrel{(2)}{F} =0\  . 
\la{1pf}\eea
Clearly, 2-partial flatness is tantamount to flatness. By Theorem \ref{pfymTh},
a 1-partially flat connection is a Yang-Mills connection.

\subsubsection{Construction of 0-partially flat connections}
Let $M$ be a $0$--admissible spin $\frac{3}{2}$ Grassmann manifold with
a 0-partially flat connection $\n$ (satisfying \re{0pf}) in a holomorphic
vector bundle $W\ra M$. The pull-back of such a connection $\n$ to a
connection in the bundle $\pi^*W\ra S_H$, where $\pi: S_H\ra M$, 
has curvature $F$ with components given by:
\bea 
F(X^e_{\pm\pm\pm} , X^{e'}_{\pm\pm\pm}) &=& 0
\nonumber\\[6pt]
F(X^{e'}_{\pm\pm\pm} , X^e_\pm) &=&
u^{\a_1}_\pm u^{\a_2}_\pm u^{\a_3}_\pm u^{\b_1}_\pm u^{\b_2}_\pm u^{\b_3}_\mp
          F(X^e_{\a_1\a_2\a_3}, X^{e'}_{\b_1\b_2\b_3})
\nonumber\\[4pt]
&=&
\pm 12\Fone_{\a_1\a_2\b_1\b_2}u^{\a_1}_{\pm}u^{\a_2}_{\pm}u^{\b_1}_{\pm}
u^{\b_2}_{\pm} =: \pm 12 \Fone_{\pm\pm\pm\pm}
\nonumber \\[6pt]
F(X^e_\pm , X^{e'}_\pm) &=&
u^{\a_1}_\pm u^{\a_2}_\pm u^{\a_3}_\mp u^{\b_1}_\pm u^{\b_2}_\pm u^{\b_3}_\mp
          F(X^e_{\a_1\a_2\a_3}, X^{e'}_{\b_1\b_2\b_3})
\nonumber\\[4pt]
&=&
-8\Ftwo_{\a_1\b_1}u^{\a_1}_{\pm}u^{\b_1}_{\pm}
=: -8\Ftwo_{\pm\pm}
\nonumber\\[6pt]
F(X^e_{+++} , X^{e'}_{---}) &=&
u^{\a_1}_+ u^{\a_2}_+ u^{\a_3}_+ u^{\b_1}_- u^{\b_2}_- u^{\b_3}_-
          F(X^e_{\a_1\a_2\a_3}, X^{e'}_{\b_1\b_2\b_3})
\nonumber\\[4pt]
&=& 36 \left( \Fone_{\a_1\a_2\b_1\b_2}u^{\a_1}_+ u^{\a_2}_+ u^{\b_1}_-
u^{\b_2}_- +  \Ftwo_{\a_1\b_1}u^{\a_1}_+u^{\b_1}_- + \Fthree\right)
\nonumber\\[4pt]
&=:& 36 \left(\Fone_0 +\Ftwo_0 +\Fthree_0 \right)
\nonumber\\[6pt]
F( X^e_{\pm\pm\pm} , X^{e'}_{\mp}) &=&
u^{\a_1}_\pm u^{\a_2}_\pm u^{\a_3}_\pm u^{\b_1}_\pm u^{\b_2}_\mp u^{\b_3}_\mp
          F(X^e_{\a_1\a_2\a_3}, X^{e'}_{\b_1\b_2\b_3})
\nonumber\\[4pt]
&=& \pm 24 \Fone_{\pm \pm} + 12 \Ftwo_{\pm \pm}
\nonumber\\[6pt]
F(X^e_+ , X^{e'}_-) &=&
u^{\a_1}_+ u^{\a_2}_+ u^{\a_3}_- u^{\b_1}_+ u^{\b_2}_- u^{\b_3}_-
          F(X^e_{\a_1\a_2\a_3}, X^{e'}_{\b_1\b_2\b_3})
\nonumber\\[4pt]
&=&  12\Fone_0 -4 \Ftwo_0 - 12 \Fthree_0
\nonumber\\[4pt]
F(v\,, \,\cdot\,  ) &=& 0\ ,
\label{componentsEqu}\eea
where $v,v'$ are vertical vector fields on $S_H$. The form of these
components lead us to:
\bd
A connection in a holomorphic vector bundle $W\ra S_H$ is
{\bss 0-partially flat} if its curvature satisfies the equations,
\bean 
F(X^e_{\pm\pm\pm} , X^{e'}_{\pm\pm\pm})&=& 0
\\[6pt]
F(X^{e}_{\pm\pm\pm} , X^{e'}_{\pm})&=& F(X^{e'}_{\pm\pm\pm} , X^{e}_{\pm})
\\[6pt]
F(X^{e}_{\pm} , X^{e'}_{\pm})&=& - F(X^{e'}_{\pm} , X^{e}_{\pm})
\\[6pt]
F(v\,, \,\cdot\,  )&=&  0\quad ,\quad\forall v\in T^vS_H\,.
\eean
\ed
The restriction of a 0-partially flat connection to a leaf of the integrable
distribution ${\langle }\cd_{\!3+}, \p_0{\rangle }$ is clearly flat. 
In this case, an {\bss analytic frame} in the holomorphic vector
bundle $\pi^*\nu: \bC^r \times S_H \ra S_H$ is a frame which is parallel
along leaves of this integrable distribution. With respect to such a frame, 
a connection in the vector bundle $\pi^*\nu$ can be written as
$$\left\{\arr 
&\nz &=\  \p_0  \\[6pt]
&\n^S_{X^e_{+++}} &=\  X^e_{+++} \\[6pt]
&\n^S_{X^e_\pm} &=\  X^e_\pm + A(X^e_\pm)  \\[6pt]
&\npm &=\   \dpm + A(\dpm ) \\[6pt]
&\n^S_{X^e_{---}} &=\  X^e_{---} +  A(X^e_{---})\ .
\ea\right.$$
\bd
A connection $\n^S$ over $S_H$ is called {\bss almost 0-partially flat} if
its curvature satisfies the following equations:
\bea
&& F(X^e_{+++}, X^{e'}_{+++}) = F(X^{e}_{+++},v) 
=  F(X^{e}_\pm ,v)= 0 \quad,\quad \forall v\in T^vS_H\ ,
\nonumber\\[4pt]
&& F(\dpp,\,\cdot\,  ) = F(\p_0,\,\cdot\, )=  0\, . 
\label{almost0Equ}\eea
\ed

Following the construction of almost half-flat connections, we may
construct almost 0-partially flat connections, which allow
deformation to a 0-partially flat connection. 
As in the case of a half-flat connection (c.f.\ Proposition \ref{ahfProp1}),
an almost 0-partially flat connection is completely determined by
the potentials  $A_{\pm\pm} =: A(\dpm )$ with respect to an analytic frame. 
\bp \label{a0pfProp1}
Let $\n$ be an almost 0-partially flat connection in the vector bundle
$\pi^*\nu: \bC^r {\times} S_H \ra S_H$  with potentials
$A_{++}$, $A_{--}$, $A(X_\pm^e)$ and $A(X_{---}^e)$ in an analytic frame. 
Then:\\
\noindent
{\it (i)\ } The potential $A_{++}$ is analytic and has charge $2$, i.e.
\be \label{0pfprepEqu}
X^e_{+++} \App = 0\quad ,\quad  \p_0 \App = 2 \App\ . 
\ee
{\it (ii)\ }  The potential $A_{--}$ satisfies
\be
\dpp \Amm - \dmm \App + [ \App , \Amm ] = 0 \quad,\quad
\p_0 \Amm =-2 \Amm\ . 
\la{0pfamm}\ee
{\it (iii)\ } The potentials $A(X_\pm^e)$  and $A(X_{---}^e)$ are
then recursively determined as follows:
\bea
   A(X^e_+) &=&  - \fr13 X^e_{+++} \Amm
\nonumber\\
   A(X_-^e) &=&
  \fr12 \left(\dmm A(X_+^e) -X_+^e\Amm + [\Amm , A(X_+^e)] \right)
\nonumber\\
A(X_{---}^e) &=&
\dmm A(X_-^e) -X_-^e\Amm + [\Amm , A(X_-^e)]
\label{0pfpots} \eea
and they have charges $+1,-1$ and $-3$ respectively, i.e. 
\be
\p_0 A(X^e_\pm) =\pm A(X^e_\pm) \quad ,\quad
\p_0 A(X_{---}^e) =-3 A(X_{---}^e)\ .
\la{0pfch}\ee
Conversely, any set of matrix-valued potentials $A_{++}$, $A_{--}$, $A(X_\pm^e)$
and $A(X_{---}^e)$
satisfying \re{0pfprepEqu}--\re{0pfch} define an
almost 0-partially flat connection. 
\ep

\pf
(i) The curvature constraints
$F(X^e_+\,,\,\dpp)=0\ ,\ F(\p_0\,,\,\dpp)=0$, in an analytic frame,
take the form \re{0pfprepEqu}. \\
(ii) The further almost 0-partial-flatness conditions,
$F(\dpp\,,\,\dmm) = F(\p_0\,,\,\dmm)=0$ give
equations \re{0pfamm} for the potential $\Amm\,$. \\
(iii) Having obtained $\,\Amm\,$, we can find
$A(X_\pm^e)$  and $A(X_{---}^e)$ from
the equations
\bean
F(\dmm , X_{+++}^e) = 0 &\Leftrightarrow &
-X_{+++}^eA_{--} =  A([\dmm , X_{+++}^e]) = 3A(X_+^e) \\[6pt]
F(\dmm , X_+^e) = 0  &\Leftrightarrow &
\dmm A(X_+^e) {-} X_+^e\Amm {+} [\Amm , A(X_+^e)] {=} A([\dmm , X_+^e]) =
2 A(X_-^e)\\[6pt]
F(\dmm , X_-^e) = 0  &\Leftrightarrow &
\dmm A(X_-^e) {-} X_-^e\Amm {+} [\Amm , A(X_-^e)] {=} A([\dmm , X_-^e]) =
A(X_{---}^e) \,. 
\eean
Equations \re{0pfch} follow from \re{0pfpots}. 
\qed

\noindent
Now, using this proposition, a modification of Theorem \ref{algo} gives an
algorithm for the construction of all almost 0-partially flat connections. 

\bt Let $\App$ be an analytic prepotential, i.e.\  a matrix-valued function
on a domain $U = \pi^{-1}(V) \subset S_H$, $V\subset M$ a simply
connected domain, satisfying \re{0pfprepEqu}, and $\F$ an
invertible matrix-valued function on $U$ which satisfies the
equations
\be \label{0pfb} \dpp \Phi = -A_{++}\Phi \, ,\quad \p_0 \Phi = 0\, . 
\ee
Such a function $\F$ always exists. 
Then the pair $(\App , \F)$ determines an almost 0-partially flat connection
$\n^S = \n^{(\App , \F)}$. Its potentials with respect to an analytic frame
are given by $\App$, $\Amm = -(\dmm \F) \F^{-1}$ and \re{0pfpots}. 
Conversely, any almost  0-partially flat connection is of this form. 
\et
The proof follows that for Theorem \ref{algo} and uses
Proposition \ref{a0pfProp1}. 

To deform an almost 0-partially flat connection into a
0-partially flat connection, we need to find a transformation from
the above analytic frame to a central frame. 
Analogously to Lemma \ref{centralLem} we may prove:

\bl\la{0pfcentralLem}
Let $\n = \n^{(A_{++},\F)}$ be the almost 0-partially flat
connection associated to the analytic prepotential $A_{++}$ with respect to
the analytic frame $\f$ and an invertible solution $\F$ of $\re{0pfb}$. 
Then the frame $\psi := \f\F$ is a central frame for the connection $\n$,
i.e.\  the potentials $C(\dpm )$ and $C(\p_0)$ with respect to that frame
vanish. 
\el
With respect to the central frame $\psi$, the above almost 0-partially flat
connection then takes the form:
$$\left\{ \arr
&\n^S_{X^e_{+++}} &=\ X^e_{+++} + C({X^e_{+++}})\
                  =\  X^e_{+++} + \F\ X^e_{+++}\ \F^{-1}
\la{0pfp}\\[6pt]
&\n^S_{X^e_\pm} &=\ X^e_\pm + C({X^e_\pm})\\[6pt]
&\n^S_{X^e_{---}}&=\  X^e_{---} + C({X^e_{---}})  \\[6pt]
&\npp &=\   \dpp  \quad,\quad \nmm\ =\   \dmm \quad,\quad 
\nz\ =\  \p_0 \, . 
\ea \right. $$
where in terms of the analytic frame potentials $A(X)$, the central frame
potentials $C(X)$ are given by
$ C(X) = \F^{-1}A(X)\F + \F^{-1}(X\F)$. 
Moreover, the equations $F(\dpp\,,\,X^e_{+++})= F(\p_0\,,\,X^e_{+++})=0$ 
imply that the potential $C({X^e_{+++}})$ satisfies the equations
\be
\dpp C({X^e_{+++}}) = 0\ ,\quad \p_0 C({X^e_{+++}}) = 3 C({X^e_{+++}})\, . 
\la{0pfCp}\ee
The following proposition is analogous to Proposition \ref{linC} in the
half-flat case. 
\bp\la{3C}
The potential $C({X^e_{+++}})$ of an almost 0-partially flat connection $\n$
with respect to a central frame is cubic in $u_+^\a\,$,
\be
C({X^e_{+++}}) = u_+^\a u_+^\b u_+^\g  C(\wt{X^e_{\a\b\g}})
               = u_+^\a u_+^\b u_+^\g  C^e_{\a\b\g} \ , 
\ee
where the coefficients $C^e_{\a\b\g} = C^e_{\a\b\g}(x^i)$, symmetric in
$\a,\b,\g\,$, are matrix valued functions of coordinates $x^i$ on $M$
and $(x^i, u_\pm^\a)$ are the local coordinates associated with the
trivialisation $S_H = M\times \Sp{1,\bC}$.
\ep

\noindent
With respect to a central frame, we can therefore write
$\n_{X^e_{+++}} = X^e_{+++} + u^\a_+ u_+^\b u_+^\g  C^e_{\a\b\g}\,$.  
Using $C^e_{\a\b\g}$, we now define a new connection in $\pi^* \nu$ 
over $S_H$ by 
$$ \left\{ \arr
&\wh\n_{X^e_{+++}} &=\ X^e_{+++} + u^\a_+ u_+^\b u_+^\g  C^e_{\a\b\g}  \\[6pt]
&\wh\n_{X^e_+} &=\  X^e_+ + u^\a_+  u_+^\b u_-^\g  C^e_{\a\b\g}  \\[6pt]
&\wh\n_{X^e_-} &=\  X^e_- + u^\a_- u_-^\b u_+^\g  C^e_{\a\b\g}   \\[6pt]
&\wh\n_{X^e_{---}} &=\  X^e_{---} + u^\a_- u_-^\b u_-^\g  C^e_{\a\b\g} \\[6pt]
&\wh\n_{\dpp} &=\  \dpp \quad,\quad \wh\n_{\dmm}\ =\ \dmm \quad,\quad
\wh\n_{\p_0}\ =\ \p_0 \, . 
\ea \right. $$

\noindent
The following theorem is the analogue of  Theorem \ref{hatThm} in the
half-flat case. 
\bt\label{0pfcentral}
The constructed connection $\wh\n$ is a 0-partially flat connection in
$\pi^*\nu$ over $S_H$
and it is the pull-back of the following 0-partially flat connection
$\n^M$ in $\nu$ over $M$:
$$
\n^M_{X^e_{\a\b\g}} = X^e_{\a\b\g}  + C^e_{\a\b\g}\ . 
$$
\et

\pf
As in Lemma \ref{pb} we may show that the connection $\wh\n$ is the pull-back
of the connection $\n^M$. It then suffices to show that $\n^M$ is
0-partially flat. The connections $\n$ and $\wh\n$ coincide in the direction 
of $X^e_{+++}\,$. Hence, using 
$ C^e_{+++} := u^{\a_1}_+u^{\a_2}_+u^{\a_3}_+ C^e_{\a_1\a_2\a_3}\,$,  we have
\bean
0 &=&  F^{\n}(X^e_{+++} , X^{e'}_{+++})  \\[4pt]
     &=& F^{\wh\n}(X^e_{+++} , X^{e'}_{+++}) \\[4pt]
     &=&  X^e_{+++} C^{e'}_{+++} - X^{e'}_{+++} C^e_{+++} +
\bigl[ C^e_{+++} \,,\, C^{e'}_{+++} \bigr]
      - C\bigl({\bigl[ X^e_{+++} \,,\, X^{e'}_{+++}\bigr]}\bigr) 
\\[4pt]       
&=& u^{\a_1}_+u^{\a_2}_+u^{\a_3}_+ u^{\b_1}_+u^{\b_2}_+u^{\b_3}_+ 
\left(
X^e_{\a_1\a_2\a_3} C^{e'}_{\b_1\b_2\b_3} 
- X^{e'}_{\b_1\b_2\b_3} C^e_{\a_1\a_2\a_3} 
+ \bigl[ C^e_{\a_1\a_2\a_3} \,,\, C^{e'}_{\b_1\b_2\b_3} \bigr] 
\right.\\  && \left. \hskip 4.0truecm                    
- C\bigl({\bigl[X^e_{\a_1\a_2\a_3} \,,\, X^{e'}_{\b_1\b_2\b_3} \bigr]}\bigr) 
\right) \\[4pt]
      &=& u^{\a_1}_+u^{\a_2}_+u^{\a_3}_+u^{\b_1}_+u^{\b_2}_+u^{\b_3}_+\  
F^{\n^M}\bigl( X^e_{\a_1\a_2\a_3}\, ,\,X^{e'}_{\b_1\b_2\b_3}  \bigr)\ ,
\eean
since $X^e_{+++} u^\b_+ = 0$. 
This shows that the component $\stackrel{(0)}{F}=0$ in the decomposition
\re{3curv_decomp}, i.e.\  the connection $\n^M$ is 0-partially flat. 
\qed

\subsubsection{Construction of 1-partially flat connections}
Let $M$ be a $1$--admissible spin $\frac{3}{2}$ Grassmann manifold with
a 1-partially flat connection $\n$ (satisfying \re{1pf}) in a holomorphic
vector bundle $W\ra M$. The pull-back of such a connection to a connection
in the bundle $\pi^* W\ra S_H$, where $\pi: S_H\ra M$, has curvature $F$ 
with components given by: 
\bea
F(X^e_{\pm\pm\pm} , X^{e'}_{\pm\pm\pm}) &=& 0 \nonumber\\[4pt]
F(X^{e'}_{\pm\pm\pm} , X^e_\pm) &=& 0 \nonumber\\[4pt]
F(X^e_\pm , X^{e'}_\pm) &=& 0 \nonumber\\[4pt]
F(X^e_{+++} , X^{e'}_{---}) &=& 36 \Fthree \nonumber\\[4pt]
F(X^e_{\pm\pm\pm} , X^{e'}_{\mp}) &=& 0 \nonumber\\[4pt]
F(X^e_+ , X^{e'}_-) &=& -12 \Fthree \nonumber\\[4pt]
F(v\,, \,\cdot\, ) &=& 0\ ,
\la{1f}\eea
where $v$ is any vertical vector field on $S_H$. A connection in a
holomorphic vector bundle $W\ra S_H$ is
{\bss 1-partially flat} if its curvature satisfies the above equations. 
The restriction of a 1-partially flat connection to a leaf of the integrable
distribution ${\langle}\cd_{\!3+}, \cd_{\!+}, \p_0{\rangle}$ is clearly flat. 
In this case, an {\bss analytic frame} in the holomorphic vector
bundle $\pi^*\nu: \bC^r \times S_H \ra S_H$ is a frame which is parallel
along leaves of this distribution. 
With respect to such a frame, a connection in the vector bundle
$\pi^*\nu$ can be written as
$$\left\{ \arr
&\nz &=\  \p_0  \\[6pt]
&\n^S_{X^e_{+++}} &=\  X^e_{+++}  \\[6pt]
&\n^S_{X^e_+} &=\  X^e_+  \\[6pt]
&\npp &=\   \dpp + \App\  \\[6pt]
&\nmm &=\   \dmm + \Amm\  \\[6pt]
&\n^S_{X^e_{---}} &=\  X^e_{---} + A(X^e_{---}) \\[6pt]
&\n^S_{X^e_-} &=\  X^e_- + A(X^e_-)\ ,    
\ea \right.$$
with potentials $A(X^e_{+++})=A(X^e_+) = A(\p_0)= 0$. 
We look for solutions of the system \re{1f} in this analytic gauge. 

\bd
A connection $\n^S$ over $S_H$ is called {\bss almost 1-partially flat}
if its curvature satisfies the equations:
\bea
&& F(X^e_{+++}, X^{e'}_{+++}) = F(X^e_{+++}, X^{e'}_{+}) 
= F(X^e_{+}, X^{e'}_{+}) = 0\nonumber\\[4pt]
&& F(X^{e}_{+++},v ) =  F(X^{e}_\pm ,v ) 
= F(\dpp,\,\cdot\,  ) = F(\p_0,\,\cdot\, ) =  0\quad,\quad \forall v\in T^vS_H\ . 
\label{almost1Equ}\eea
\ed
In virtue of these equations, the potentials $A_{\pm\pm} =: A(\dpm )$
determine all other potentials:
\bp \label{a1pfProp1}
Let $\n$ be an almost 1-partially flat connection in the vector bundle
$\pi^*\nu: \bC^r \times S_H \ra S_H$  with potentials
$A_{++}$, $A_{--}$, $A(X_-^e)$ and $A(X_{---}^e)$ in an analytic frame. 
Then:\\
\noindent
(i)\  The potential $A_{++}$ is analytic and has charge $2$, i.e.\ 
\be
X^e_{+++} \App = 0\quad ,\quad X^e_+ \App = 0\quad ,\quad
 \p_0 \App = 2 \App\ . 
\label{1pfprepEqu}\ee
\noindent
(ii)\  
The potential $A_{--}$ satisfies
\be
\dpp \Amm - \dmm \App + [ \App , \Amm ] = 0 \quad,\quad
\p_0 \Amm =-2 \Amm\ . 
\la{1pfamm}\ee
\noindent
(iii)\   
The potentials $A(X_-^e)$  and $A(X_{---}^e)$ are then
recursively determined as follows:
\bea
   A(X_-^e) &=&  - \fr12  X_+^e\Amm
\nonumber\\[2pt]
A(X_{---}^e) &=&
\dmm A(X_-^e) -X_-^e\Amm + [\Amm , A(X_-^e)]
\label{1pfpots} \eea
and they have charges $-1$ and $-3$ respectively; i.e. 
\be
\p_0 A(X^e_-) =\pm A(X^e_-) \quad ,\quad
\p_0 A(X_{---}^e) =-3 A(X_{---}^e)\ .
\la{1pfch}\ee
Conversely, any set of matrix-valued potentials  $A_{++}$, $A_{--}$,
$A(X_-^e)$ and $A(X_{---}^e)$
satisfying \re{1pfprepEqu}--\re{1pfch} define an
almost 1-partially flat connection. 
\ep

\pf
(i) Equations \re{1pfprepEqu} are equivalent to
$F(\dpp , X_{+++}^e) {\,=\,} F(\dpp , X_{+}^e){\,=\,}F(\p_0, \dpp ){=}0$. 

\noindent
(ii) The further almost 1-partial-flatness conditions,
$F(\dpp\,,\,\dmm) = F(\p_0\,,\,\dmm)=0$ give equations \re{1pfamm}. 

\noindent
(iii) Having obtained $\,\Amm\,$, we can find
$A(X_-^e)$  and $A(X_{---}^e)$ from the equations
\bean 
F(\dmm , X_+^e) = 0  &\Leftrightarrow &
-X_+^e\Amm  = A([\dmm , X_+^e]) = 2 A(X_-^e)\\[3pt]
F(\dmm , X_-^e) = 0  &\Leftrightarrow &
\dmm A(X_-^e) -X_-^e\Amm + [\Amm , A(X_-^e)] = A([\dmm , X_-^e]) =
A(X_{---}^e) \,. \eean
The equations \re{1pfch} follow from \re{1pfpots}. 
\qed

\noindent
Now, starting from a prepotential $A_{++}\,$, which solves \re{1pfprepEqu}, 
we may construct an almost 1-partially flat connection. 
The potential $A_{--} = - (\dmm \Phi ) \Phi^{-1}$ is
determined, as before, from a solution $\F$ of \re{0pfb}. 
Then, with the remain potentials in an analytic frame being given by
\re{1pfpots} and satisfying \re{1pfch},
all the other equations in  \re{almost1Equ} follow. This shows that a
almost 1-partially flat connection is determined
by an arbitrary analytic prepotential $A_{++}$ and an invertible solution
$\F$ of  \re{0pfb}. 
As before, $\F$ is a transition function
from an analytic frame to a central frame, in which the above
almost 1-partially flat
connection  takes the form:
$$\left\{ \arr
&\n^S_{X^e_{+++}} &=\ X^e_{+++} + C({X^e_{+++}})\
                  =\ X^e_{+++} + \F\ X^e_{+++}\ \F^{-1} \\[8pt]
&\n^S_{X^e_+} &=\ X^e_+ + C({X^e_+})\
                  =\ X^e_+ + \F\ X^e_+\ \F^{-1} \la{1pfp}\\[8pt]
&\n^S_{X^e_-}&=\  X^e_- + C({X^e_-})  \\[8pt]
&\n^S_{X^e_{---}}&=\  X^e_{---} + C({X^e_{---}})  \\[8pt]
&\npp &=\   \dpp  \quad,\quad
\nmm =   \dmm  \quad,\quad
\nz =  \p_0\ .
\ea \right. $$
Moreover, the equations $F(\dpp\,,\,X^e_{+++})=F(\p_0\,,\,X^e_{+++})=0$
imply that the potential $C({X^e_{+++}})$ satisfies the equations
$$
\dpp C({X^e_{+++}}) = 0\ ,\quad \p_0 C({X^e_{+++}}) = 3 C({X^e_{+++}})\ .
$$

\bp
The potentials $C({X^e_{+++}})$ and $C({X^e_{+}})$ of an almost
1-partially flat connection $\n$
with respect to a central frame have the form,
$$
C({X^e_{+++}}) = u_+^\a u_+^\b u_+^\g  C^e_{\a\b\g}\quad,\quad
C({X^e_{+}}) = u_+^\a  u_+^\b u_-^\g C^e_{\a\b\g}\, ,
$$
where $C^e_{\a\b\g}$ is a function on $M$, symmetric in $\a,\b,\g$. 
\ep

\noindent
With respect to a central frame, we can therefore write
$$
\n_{X^e_{+++}} = X^e_{+++} + u^\a_+ u_+^\b u_+^\g  C^e_{\a\b\g}\quad,\quad
\n_{X^e_+} = X^e_+ + u_+^\a  u_+^\b u_-^\g C^e_{\a\b\g}\, . 
$$
We define a modified connection in the bundle $\pi^* \nu$ over $S_H$ by
\be \left\{ \arr
&\wh\n_{X^e_{+++}} &=\  X^e_{+++} + u^\a_+ u_+^\b u_+^\g  C^e_{\a\b\g} 
\\[8pt] &\wh\n_{X^e_+} &=\  X^e_+ + u^\a_+  u_+^\b u_-^\g  C^e_{\a\b\g} 
\\[8pt] &\wh\n_{X^e_-} &=\  X^e_- + u^\a_- u_-^\b u_+^\g  C^e_{\a\b\g}  
\\[8pt] &\wh\n_{X^e_{---}}&=\  X^e_{---} + u^\a_- u_-^\b u_-^\g C^e_{\a\b\g}  
\\[8pt] &\wh\n_{\dpp} &=\   \dpp \quad,\quad \wh\n_{\dmm}\ =\ \dmm
\quad,\quad \wh\n_{\p_0}\ =\ \p_0 \, .  \ea \right. \ee

\noindent
As in the 0-partially flat case, we have:

\bt
The constructed connection $\wh\n$ is a 1-partially flat connection in
$\pi^*\nu$ over $S_H$
and it is the pull-back of the following 1-partially flat connection
$\n^M$ in $\nu$ over $M$:
\be
\n^M_{X^e_{\a\b\g}} = X^e_{\a\b\g}  + C^e_{\a\b\g}\, . 
\la{1pfM}\ee
\et

\pf
As before one shows that the connection $\wh\n$ is the pull-back of
the connection $\n^M$. It remains to show that $\n^M$
is 1-partially flat. Since any almost 1-partially flat connection
is almost 0-partially flat we have $\stackrel{(0)}{F}=0$ by
Theorem \ref{0pfcentral}. Next we show
that $\stackrel{(1)}{F}=0$. The connections $\n$ and $\wh\n$ coincide in
the direction of $X^e_{+++}$ and
$X_+^e$. Hence, using equation \re{componentsEqu}, which holds for
0-partially flat connections, we have
$$\arr
0 &=&  F^{\n}(X^e_{+++} , X^{e'}_{+})  \\[6pt]
     &=& F^{\wh\n}(X^e_{+++} , X^{e'}_{+}) \\[6pt]
     &=& 12 \Fone_{\a_1\a_2\b_1\b_2}u^{\a_1}_{+}u^{\a_2}_{+}u^{\b_1}_{+}
u^{\b_2}_{+}\, . 
\ea$$
This shows that the component $\stackrel{(1)}{F}$ in the decomposition
\re{3curv_decomp} vanishes. Similarly,
$$\arr
0 &=&  F^{\n}(X^e_{+} , X^{e'}_{+})  \\[6pt]
     &=& F^{\wh\n}(X^e_{+} , X^{e'}_{+}) \\[6pt]
     &=& -8\Ftwo_{\a_1\b_1}u^{\a_1}_{+}u^{\b_1}_{+}\, ,
\ea$$ implies $\Ftwo = 0$, and hence that $\wh\n$ is 1-partially
flat. \qed

\noindent
By Theorem \ref{pfymTh}, the 1-partially flat connection 
$\n^M_{X^e_{\a\b\g}}$ in \re{1pfM} is a Yang-Mills connection.

\vfil
%%%%%%%%%%%%%%%%%%%%%%%%%%%%%%%%%%%%%%%%%%%%%%%%%%%%%%%%


\begin{thebibliography}{XXXX}
\baselineskip=16pt
\addtolength{\itemsep}{-5pt}

\bibitem[AG]{AG}  D.V.\ Alekseevsky and M.M.\ Graev, 
{\it Grassmann and hyperk\"ahler structures on some spaces of sections 
of holomorphic bundles} in {\sl  Manifolds and geometry},
Ed.\ de Bartolomeis et al, CUP, Cambridge, 1996

\bibitem[AHS]{AHS} M.F.\ Atiyah, N.J.\ Hitchin and  I.M.\ Singer,
{\it Self-duality in four-dimensional Riemannian geometry},
{\sl Proc.\ Roy.\ Soc.\ Lond.\/ \bf A362} (1978) 425--461
%%CITATION = PRSLA,A362,425;%%

\bibitem[B]{B} A.\ L.\ Besse, {\it Einstein manifolds}, Springer, Berlin, 1987

\bibitem[BKS]{BKS}
L.\ Baulieu, H.\ Kanno and I.M.\ Singer,
{\it Special quantum field theories in eight and other dimensions},
{\sl Commun.\ Math.\ Phys.\/  \bf 194} (1998) 149--175
[hep-th/9704167]
%%CITATION = HEP-TH 9704167;%%


\bibitem[CS]{CS} M.M.\ Capria  and  S.M.\ Salamon,
               {\it Yang-Mills fields on quaternionic spaces},
               {\sl Nonlinearity \/ \bf 1} (1988) 517--530
%%CITATION = NOLIN,1,517;%%

\bibitem[CDFN]{CDFN} E.\ Corrigan, C.\ Devchand, D.B.\ Fairlie and J.\ Nuyts,
                {\it First order equations for gauge fields in spaces of
                 dimension greater than four},
                 {\sl Nucl.\ Phys.\/  \bf B214} (1983) 452--464
%%CITATION = NUPHA,B214,452;%% 

\bibitem[CGK]{CGK} E.\ Corrigan, P.\ Goddard and A.\ Kent, 
{\it Some comments on the ADHM construction in $4k$ dimensions},  
{\sl Commun.\ Math.\ Phys.\ \/ \bf 100} (1985) 1-13
%%CITATION = CMPHA,100,1;%%


\bibitem[DN1]{DN1} C.\ Devchand and J.\ Nuyts,
{\it Supersymmetric Lorentz-covariant hyperspaces and self-duality equations  
in dimensions greater than (4$|$4)},
{\sl Nucl.\ Phys.\/ \bf B503} (1997) 627--656
[hep-th/9704036]
%%CITATION = HEP-TH 9704036;%%

\bibitem[DN2]{DN2} C.\ Devchand and J.\ Nuyts,
{\it Super self-duality for Yang-Mills fields in dimensions greater than 
four}, {\sl J. High Energy Phys.\/ \bf 12} (2001) 020 [hep-th/0109072]
%%CITATION = HEP-TH 0109072;%%

\bibitem[DT]{DT} S.K.\ Donaldson and R.P.\ Thomas, 
{\it Gauge theory in higher dimensions}, in
{\sl The geometric universe: science, geometry,
and the work of Roger Penrose}, ed. S.A.\ Huggett, et al., 
Oxford Univ.\ Press, 31-47 (1998)

\bibitem[GIOS]{GIOS} 
A.\ Galperin, E.\ Ivanov, V.\ Ogievetsky and E.\ Sokatchev,
{\it Harmonic superspace}, Cambridge Univ.\ Press, U.K., 2001; 
{\it Gauge field geometry from complex and harmonic analyticities. 
I. K\"ahler and selfdual Yang-Mills cases},  
{\sl  Ann. Phys.\/ \bf 185} (1988) 1--21 
%%CITATION = APNYA,185,1; %%


\bibitem[K]{K} B.\ Kostant,
{\it On invariant skew-tensors},
{\sl Proc.\ Nat.\ Acad.\ Sci.\ USA.\/ \bf 42 } (1956), 148--151
%%CITATION = PNASA,42,148;%%

\bibitem[N]{N} T.\ Nitta, {\it Vector bundles over quaternionic K\"ahler
manifolds}, {\sl Tohoku Math.\ J.\/ \bf 40} (1988) 425--440 
%%CITATION = TOMJA,40,425;%%
  
\bibitem[OV]{OV} A.L.\ Onishchik  and  E.B.\ Vinberg, 
{\it Lie Groups and Algebraic Groups},  
Springer-Verlag, Berlin, 1990 

\bibitem[T]{T} G.~Tian, {\it Gauge theory and calibrated geometry. I.},
{\sl  Annals Math. (2) \/ \bf 151} (2000) 193--268 
%%CITATION = MATH-DG 0010015;%%

\bibitem[W]{W} R.~S.~Ward,
{\it Completely solvable gauge field equations in dimension greater than four},
{\sl Nucl.\ Phys.\/ \bf B236} (1984) 381--396
%%CITATION = NUPHA,B236,381;%%

\end{thebibliography}
\end{document}